%% file: RandomCommutingMatrices_13May20254.tex
\documentclass[12pt]{article}
\usepackage{graphicx}
\usepackage{multicol,multirow}
\usepackage{amsmath,amssymb,amsfonts}
\usepackage{mathrsfs}
\usepackage{rotating}
\usepackage{appendix}
\usepackage[numbers]{natbib}

\usepackage{verbatim}
\usepackage{amssymb}
\usepackage{amsthm}
\usepackage{amsmath}
\usepackage{mathrsfs}
\usepackage[margin=1in]{geometry}
\usepackage{mathtools}

\input def_latex5

\usepackage[usenames,dvipsnames] {color}

\definecolor{dark_purple}{rgb}{0.4, 0.0, 0.4}
\definecolor{dark_green}{rgb}{0.0, 0.7, 0.0}

\numberwithin{equation}{section}

\title{Random commuting matrices}
\author{
John E. M\raise.5ex\hbox{c}Carthy
\thanks{Partially supported by National Science Foundation Grant  
DMS 2054199}
}
\date{\today}

\newcommand{\mc}{M\raise.45ex\hbox{c}Carthy}

\usepackage{mathtools}

\newcommand\mn{{\mathbb M}_n}
\newcommand\mnd{{\mathbb M}_n^d}

\renewcommand\l{\lambda}
\newcommand\vdn{{\mathfrak C}^d_n}
\newcommand\vdt{{\mathfrak C}^d_2}
\newcommand\vtn{{\mathfrak C}^2_n}
\newcommand\vhdn{{\mathfrak C}^{d,H}_n}
\newcommand\vhdpn{{\mathfrak C}^{d',H}_n}
\newcommand\vhdmrn{{\mathfrak C}^{d-1,H}_{r_i}}
\newcommand\comhdng{{\mathfrak C}^{d,H}_{n,{\rm gen}}}
\newcommand\comdng{{\mathfrak C}^{d}_{n,{\rm gen}}}
\newcommand\comdn{{\mathfrak C}^d_n}
\newcommand\comtn{{\mathfrak C}^2_n}
\newcommand\comhdn{{\mathfrak C}^{d,H}_n}
\newcommand\acomdn{{\mathfrak A}^d_n}
\newcommand\acomhdn{{\mathfrak A}^{d,H}_n}
\newcommand\acomhtn{{\mathfrak A}^{2,H}_n}
\newcommand\niln{{\mathfrak N}_n}
\newcommand\tdn{{\mathcal T}^d_n}
\newcommand\ton{{\mathcal T}^1_n}
\newcommand\tdt{{\mathcal T}^d_2}

\newcommand\tw{{\tilde{w}}}
\newcommand\un{{\mathcal U}(n)}
\newcommand\cp{{\mathcal P}}

\newcommand\rdn{{({\mathbb R}^d)^n}}
\newcommand\cdn{{({\mathbb C}^d)^n}}
\newcommand\dimc{\dim_\C}

\newcommand\azn{{\mathcal A}_{0,n}}
\newcommand{\fno}{F_{n,1}}
\newcommand{\fnz}{F_{n,0}}
\newcommand\diagdn{{\rm Diag}^d_n}

\newcommand\ccn{{\beta_n}}
\newcommand\en{{\mathcal E}_n}
\newcommand\rd{{\R^d}}
\newcommand\ant{A_{n,2\eta}}
\newcommand\vn{{\mathfrak V}_n}

\begin{document}

\bibliographystyle{plain}

\maketitle
\begin{abstract}
We define a {\em  random commuting $d$-tuple of $n$-by-$n$ matrices} to be a random variable that takes values
in the set of commuting $d$-tuples and has a distribution that is a rapidly decaying continuous weight on this algebraic set.
In the Hermitian case, we characterize the eigenvalue distribution, and derive the limit as $n$ tends to infinity for a Gaussian weight.
In the non-Hermitian case, we get a formula that holds if the set is irreducible. We show that there are qualitative differences between the single matrix case and the several commuting matrices case.
\renewcommand{\thefootnote}{\fnsymbol{footnote}} 
\footnotetext{\emph{AMS Subject Class:}  60B20, 15B52.}     
\renewcommand{\thefootnote}{\arabic{footnote}} 
\end{abstract}
\section{Introduction}
\label{secin}
\subsection{Overview}

The study of random matrices, and in particular their eigenvalue distribution, is a well-established field,
dating back to the pioneering work of Wigner \cite{wi55,wi57,wi58},  Dyson \cite{dy62a,dy62b, dy62c} 
and Ginibre \cite{gi65}.
For a comprehensive treatment, see
 e.g. the books \cite{
agz10,
bl09, dei99, dg09, meh04,
ms17,
ps11,
ta12}. 
The field has  many applications in mathematics and theoretical physics.
Typically the entries of the matrix are chosen independently from each other (or for Hermitian matrices the upper triangular entries are chosen independently) according to some probabilitly distribution, and then one tries to characterize the distribution of specific properties of the matrix, such as the eigenvalues, or the norm, or the singular values, both in the finite case and in the limit as the size of the matrix goes to infinity. The best understood case in random matrix theory is when the entries of the matrix are chosen i.i.d. from a Gaussian distribution.

What would it mean to choose random matrices subject to some algebraic constraint?
For a specific example, what would it mean to choose a random pair $(X,Y)$ of commuting $n$-by-$n$ matrices?
One approach would be to first choose $X$ by some random method, and then to choose something in its commutant.
Generically, $X$ will have $n$ different eigenvalues, and $Y$ will then have the same eigenvectors, so it would only be necessary to describe the distribution of the eigenvalues of $Y$. But what should that distribution be? Should it depend on the eigenvalues of $X$? And should the matrices $X$ be chosen in such a way that ones with a `larger' commutant are more heavily weighted?
If the possible entries of $X$ and $Y$ came from a large but finite set, one could in principle choose a pair independently, and then
keep it if and only if the pair commutes. But if the distribution is continuous, then the probability that any given pair would commute is zero, and one does not want to condition on a zero-probability event.

The purpose of this note is first to describe a systematic way to make sense of the notion of a random $d$-tuple of matrices that
satisfies certain algebraic relations, and which can reduce to the i.i.d. Gaussian model if the relations are absent.
Secondly, we shall explore this theory in the particular case of random commuting $d$-tuples. In the Hermitian case, we can derive the limiting distribution for a Gaussian weight; in the non-Hermitian case our results are far less complete.

\subsection{What is a random $d$-tuple of matrices that satisfy algebraic relations?}
\label{subsecaab}

Let $\mn$ denote the algebra of $n$-by-$n$ complex matrices, and $\mnd$ the 
set of $d$-tuples $X = (X^1, \dots , X^d)$ of $n$-by-$n$  matrices.
We shall use superscripts to denote the components of $X$,  subscripts to denote  sequences,
and paired subscripts the matrix entries (so $X^r_{ij}$ is the $(i,j)$ element of the $r^{\rm th}$ matrix in the $d$-tuple $X$).

Let ${\frak V}_n \subseteq \mnd$ be a set that is unitarily invariant, by which we mean that if $U$ is in the $n$-by-$n$ unitary group, denoted $\un$, and we define $UX U^*$ componentwise:
\[
UXU^* \ :=\ ( U X^1 U^*, \dots, U X^d U^*)  ,\]
then if $X \in {\frak V}_n$, so is $UXU^*$.
Identifying $\mnd$ componentwise with $\R^{2n^2d}$, we will always assume that $\vn$
is locally a smooth submanifold of $\mnd$ except on a singular set of lower dimension.

One collection of examples is given by non-commutative 
 algebraic sets, by which we mean there are non-commutative polynomials
$p_1, \dots, p_N$ in the $2d$ variables  $\{x^1, (x^1)^*, \dots, x^d, (x^d)^*\}$ so that
\be
\label{eqa01}
{\frak V}_n \= \{ X \in \mnd : p_j(X) = 0 \ \forall\ 1 \leq j \leq N \} .
\ee
(Note that the polynomials $p_j$ are applied to the matrices, not the
entries of the matrices.) Any set ${\frak V}_n$ defined as in \eqref{eqa01} is unitarily invariant.
We could also look at other conditions that are unitarily invariant, like having determinant one, or trace zero.

\begin{example}
\label{exa1}
{\rm
Our first example is the set of commuting $d$-tuples, which we shall call $\comdn$.
\[
\comdn \= \{ X \in \mnd: X^r X^s = X^s X^r \ \forall \ 1 \leq r, s \leq d\} .
\]
The commuting Hermitian $d$-tuples we shall call $\comhdn$:
\[
\comhdn \=
\{ X \in \comdn : (X^r)^* = X^r \ \forall\ 1 \leq r \leq d \} .
\]
}
\end{example}

\begin{example}
\label{exa2}
{\rm
A second example is the set of anti-commuting $d$-tuples:
\beq
\acomdn &\= &  \{ X \in \mnd: X^r X^s + X^s X^r = 0 \ \forall \ 1 \leq r \neq  s \leq d\}\\
\acomhdn &=& \{ X \in \acomdn : (X^r)^* = X^r \ \forall\ 1 \leq r \leq d \} .
\eeq
}
\end{example}

\begin{example}
\label{exa3}
 {\rm
 We can take $d = 1$, and consider the set of nilpotent matrices:
\[
\niln \= \{ X \in \mn : (X)^{n} = 0 \} .
\]
}
\end{example}

We want to define a probability measure on ${\frak V}_n$, and we want the measure to be unitarily invariant, since we do not want to privilege any particular choice of coordinates.
Let the real dimension of ${\frak V}_n$, which we will think of as a subset of $\C^{dn^2}$ with the usual Euclidean metric, be $d_{{\frak V}_n}$.
Let $dX |_{{\frak V}_n}$, which we shall just write as $dX$ when ${\frak V}_n$ is understood, be Hausdorff measure of dimension $d_{{\frak V}_n}$ on ${\frak V}_n$. 
Let $w_n$ be a unitarily invariant weight on ${\frak V}_n$, by which we mean that $w_n : {\frak V}_n \to [0,\i)$ is a non-negative measurable function that satisfies 
\beq
w_n(UXU^*) &\=& w_n(X) \qquad \forall \ U \in \un \\
\int_{{\frak V}_n} w_n(X) dX &\=& 1 .
\eeq
Let $\mu_n$ be the probability measure $w_n(X) dX$ on $V$.

\begin{definition}
\label{defa1}
A random $d$-tuple of  $n$-by-$n$ matrices in ${\frak V}_n$ is a random variable with values in ${\frak V}_n$ and distribution $\mu_n$. 
\end{definition}
The most important example for us is when the weight is Gaussian, which means it is of the form
\be
\label{eqa02}
w_n(X) \= c_n e^{-\gamma {\rm tr}[ \sum_{r=1}^d (X^r)^* X^r]} \=
c_n e^{- \gamma \sum_{r=1}^d \sum_{i,j=1}^n | X^r_{ij} |^2} ,
\ee
where $c_n$ is a normalizing constant.

\begin{example}
{\rm
If $d = 1$ and we take  ${\frak V}_n$ to be  the Hermitian matrices $\mn^H$, then with the weight
\eqref{eqa02} we will recover the Gaussian Unitary Ensemble. With  ${\frak V}_n = \mn$ we would get the
complex Gaussian matrix ensemble. With $\vn = \un$, we get Haar measure.
}\end{example}

\subsection{Why are random tuples of matrices interesting?}

When $d=2$, the study of random pairs of commuting Hermitian matrices is the same as the study of random normal matrices. These have been studied extensively---see e.g. \cite{ahm11, ahm15, cz98, ef05, wz00, mpt08, akm19,hw21}---both for their mathematical interest and their applications to physics.
To the best of my knowledge, the eigenvalue distribution for random commuting  $d$-tuples 
has not been studied in the Hermitian case for $d>2$ or the general case for $d > 1$. 
So I shall give some reasons I think the subject is worth pursuing, and list some open questions for the Gaussian
case (\eqref{eqa02} with the normalization $\gamma = \frac 12$).

If $X = (X^1, \dots, X^d)$ is a commuting $d$-tuple of $n$-by-$n$ matrices, then an eigenvector of $X$ is a non-zero vector
$v \in \C^n$ such that $X^j v = \lambda^j v$ for each $1 \leq j \leq d$. We call the $d$-tuple $ (\l^1, \dots, \l^d)$
an eigenvalue of $X$; it is a point in $\C^n$ (or $\R^n$ if each $X^j$ is self-adjoint).  If $X$ has $n$ distinct eigenvalues,
$\l_1, \dots, \l_n$, then the collection of all eigenvalues is a point in $(\R^d)^n$.

\begin{enumerate}
\item
Large random Hermitian matrices were studied by Wigner and Dyson as models for the Hamiltonian of a complicated system with many weakly interacting parts, such as heavy atoms, and the eigenvalues were the possible energy levels. 
If one wished to study a collection of $d$ such objects whose Hamiltonians commuted
then tuples of random Hermitian matrices would provide a model.

\item
Compact Lie groups have Haar measure as a natural probability measure. The construction in Subsection \ref{subsecaab}
allows one to define a probability measure on Lie algebras that are unitarily invariant. It is not so canonical, but at least for the Lie algebra of the unitary group (which is the skew-Hermitian matrices) one does get the Gaussian Unitary Ensemble (up to a factor of $i$).

\item
If one thinks of random $d$-tuples of commuting Hermitian matrices as some sort of `bosonic' model,
then one can think of random $d$-tuples of anti-commuting Hermitians as a `fermionic' model.
For the case of a Hermitian anti-commuting pair $(X,Y)$, the eigenvalues come naturally as
pairs $(\pm \lambda, \pm \mu)$ where $\lambda$ and $\mu$ are non-negative eigenvalues of $X$ and $Y$ respectively,
so can be identified with points in $(\R_+)^2$. They are not joint eigenvalues; there is a two-dimensional subspace on which $X$ and $Y$ look like
\[
\begin{pmatrix}
\lambda & 0 \\
0 & - \lambda
\end{pmatrix},
\quad
\begin{pmatrix}
0 & \mu \\
\mu& 0 
\end{pmatrix} .
\]
For finite $n$, the distribution was studied in \cite{mccmcc24}, and the pairs $(\lambda, \mu)$ were called the skew spectrum.

\begin{question}
\label{bq2}
What is the limiting distribution of the skew spectrum in $\acomhtn$? Is it area measure on a quarter disk?
\end{question}

Corollary \ref{cord1} says  that if $X$ is chosen randomly in $\comhdn$ and one knows the distribution of eigenvalues for one matrix in the $d$-tuple, $X^1$ say,
 then one
can deduce the value of $d$, provided $n$ is large enough. One could imagine using this to count the number of elements in
the set by studying one randomly chosen element very carefully. Is there a similar phenomenon in the anti-commuting case?

\begin{question}
\label{bq7}
Let $X_n$ be chosen randomly in $\acomhdn$. Does the limiting distribution of eigenvalues of $X^1_n$ depend on $d$?
Can it be explicitly described?
\end{question}

\item
The non-Hermitian case is, not surprisingly, much harder. When $d > 2$, it is no longer true that
a generic\footnote{We use generic to mean a condition on elements of a set that only fails to hold on a  
subset of strictly smaller dimension.}
$d$-tuple of commuting $n$-by-$n$ matrices will have $n$ linearly independent eigenvectors \cite{gu92},
though this is true for $d=2$ \cite{mt55}.
So if one wants to study the joint eigenvalues it would make sense to restrict one's attention to the component
of $\comdn$ that do have $n$ linearly independent eigenvectors. Some partial progress is in Section \ref{secd} below,
but it does not address the limit.
\begin{question}
\label{bq3}
 What is the limiting distribution of the joint eigenvalues of $(\frac{1}{\sqrt{n}} X, \frac{1}{\sqrt{n}}  Y)$
 when $(X,Y)$ is chosen randomly in $\comtn$?
 \end{question}
The limiting distribution of the eigenvalues of a single random matrix is area measure on a disk \cite{meh04};
this is also the limiting distribution for a pair of commuting Hermitian matrices \cite{ef05}.
Is it possible that a pair of commuting matrices behaves like a $4$-tuple of commuting Hermitian matrices?
Theorem \ref{thmd1} shows that this is not the case when $d=2$, but it might still be true in the limit.
\begin{question}
\label{bq4}
Is the answer to Question \ref{bq3} area measure on a sphere in $\C^2$?
\end{question}

\item
Multivariable operator theory studies $d$-tuples of commuting operators on separable Hilbert spaces.
For many questions, it is sufficient to understand the finite dimensional case, if one can get a bound on the constants that is independent of the dimension. The random case may be more accessible and give new insights.
 For example, it is an open question whether And\^o's inequality holds with a constant in $3$ variables, i.e. whether there is a constant $C$ so that for any
$n$ the inequality 
\be
\label{eqaa1}
\sup \{ \| p(X) \| : X \in {\mathfrak C}^3_n, \max_{1 \leq r \leq 3}\| X^r \| \leq 1 \} \ \leq \ C \| p \|_{\D^3}
\ee
holds for every polynomial $p$ in 3 variables. The fact that $C$ cannot be chosen to be $1$ was proved in \cite{var74} and \cite{cradav}; 
\begin{question}
Does \eqref{eqaa1} hold most of the time, in a probabilistic sense?
\end{question}

\item
There are deep ties between random matrices, free probability, and the von Neumann algebras of free groups---see e.g.
\cite{ms17}. If $G$ is a group  that is not free, one could hope that the group von Neumann algebra could be connected to the theory of random tuples that satisfy appropriate algebraic relations.

To frame a specific question, suppose $G$ is a discrete group with $d$ generators $g_1, \dots, g_d$ and 
$N$ relations. For each relation $R_j$ of the form 
$ g_{i_1}^{\vare_1} \dots g_{i_k}^{\vare_k} = e$, where each $\vare_i = \pm 1$,
let $p_j$ be the non-commutative polynomial
$p_j(x) = x_{i_1}^{[\vare_1]} \dots x_{i_k}^{[\vare_k]} - 1$, where we write $x_i^{[\vare]}$ to mean
$x_i$ if $\vare = 1$ and $x_i^*$ if $\vare = -1$.
Let ${\mathfrak V}_n$ be the set of $d$-tuples of $n$-by-$n$ unitaries $U$ satisfying $p_j(U) = 0$ for $1 \leq j \leq N$.
On the group von Neumann algebra of $G$, let $\tau$ be the trace and let $L_{g}$ be the image of $g$ under the left regular representation on $\ell^2(G)$.

\begin{question}
\label{bq9}
Let $U_n$ be a random tuple in ${\mathfrak V}_n$. For any word of the form
$ g_{i_1}^{\vare_1} \dots g_{i_k}^{\vare_k}$, does 
$\frac{1}{n} {\rm tr}[( U_n^{i_1} )^{[\vare_1]} \dots (U_n^{i_k})^{[\vare_k]})]$
converge almost surely to $\tau [L_{g_{i_1}}^{[\vare_1]} \dots (L_{g_{i_k}})^{[\vare_k]}]$ ?
\end{question} 

More generally, does the theory of  random tuples connect to a theory analogous to free probability, but with some extra algebraic constraints, and also to the theory of von Neumann algebras of non-free groups?

\item
Some algebraic relations cannot be attained on finite dimensional matrix algebras.
For example, if $U,V \in {\mathcal U}(n)$ and $UV = e^{2\pi i \theta} VU$, then taking determinants we see
that $e^{2 \pi i n\theta} = 1$, and in particular $\theta$ must be rational. However, if $\theta$ is irrational, Pimsner and Voiculescu \cite{pv80} showed
that if $p_n/q_n$ is the sequence of best rational approximants to $\theta$, then the irrational rotation C*-algebra $A_\theta$ can be embedded into an AF algebra that is the inductive limit of ${\mathbb M}_{q_n} \oplus {\mathbb M}_{q_{n-1}}$ under 
suitable connecting maps, and there are unitaries $U_n,V_n$ on  ${\mathbb M}_{q_n}$ satisfying
$U_n V_n = e^{2\pi i p_n/q_n} V_nU_n$. (For more background and definitions of the terms used, see \cite{da96}).
Let ${\mathfrak V}_n = \{ (U,V) \in {\mathbb M}_{q_n} : U V = e^{2\pi i p_n/q_n} VU \}$.
\begin{question}
\label{bq10}
Does sampling random pairs in ${\mathfrak V}_n$ give information about $A_\theta$?
\end{question}

\item
Our proposed framework can be used even when $d=1$ to study single matrices with algebraic constraints.
In \cite{dh01}, Dykema and Haagerup studied strictly upper triangular matrices, where 
all the elements were chosen to be i.i.d. ${ N}(0, \frac{1}{\sqrt{n}})$.
They showed that the limit converges in *-moments to a circular operator.
Is this true if one loooks at all nilpotent matrices?
\begin{question}
\label{bq8}
Let $X_n$ be chosen randomly in $\niln$. Does $X_n$ converge in *-moments to a circular operator?
\end{question}
Every nilpotent matrix is unitarily equivalent to an upper triangular matrix, so choosing $X_n$ is the same as choosing a unitary
$U_n$ and an upper-triangular $T_n$ so that $X_n = U_n^* T_n U_n$. (This representation is not unique, but that doesn't matter for this discussion). The unitaries will be equally distributed with respect to Haar measure, since all our measures are
unitarily invariant. But if $X_n$ is chosen as in Definition \ref{defa1} with a Gaussian weight on $\niln$,
 the distribution of $T_n = U_n^* X_n U_n$ will be different from the one in which all the entries are i.i.d. normal.
 
\item
There are numerous quantities whose distribution one could seek associated with a $d$-tuple of matrices in addition to the eigenvalues, such as the joint spectral radius, the joint numerical range, various notions of singular numbers, etc.
Let us just focus on the norm $\| X \| := \max_{1 \leq r \leq d} \| X^r \|$.

\begin{question}
\label{bq5}
If $X_n$ is chosen randomly in $\comhdn$ does $\frac{1}{\sqrt{n}} \| X_n\|$ tend almost surely to $R_d$, where $R_d$ 
is as in Theorem \ref{thmb1}?
If so, is there some version of the Tracy-Widom law describing how quickly it converges?
\end{question}

The numbers $R_d$ decrease from $2$ to $1$ as $d$ increases from $1$ to $4$, but then stay constant. Do the norms of
random commuting non-Hermitian matrices also tend to a limit that is the same for all $d$ above some threshold?

\begin{question}
\label{bq6}
If $X_n$ is chosen randomly in $\comdn$ does $\frac{1}{\sqrt{n}} \| X_n\|$ tend almost surely to some number $T_d$?
Is $T_d$ constant for $d \geq d_0$?
\end{question}

\item
The rate at which the eigenvalue distribution of a random $n$-by-$n$ matrix converges to the limit has been
extensively studied---see e.g. \cite{jo98} for the case of a random Hermitian matrix, and \cite{ahm11, ahm15} for a random normal matrix. In both cases the fluctuation from the limit is shown to be a Gaussian variable.
When $d \geq 4$, Theorem \ref{thmc5} says that the eigenvalue distribution converges to surface area measure on
the sphere $S^{d-1}$.
On the interior of the ball, since the limiting distribution is $0$, the fluctuation must be non-negative.
\begin{question}
\label{bq1}
At what rate does the  distribution of joint eigenvalues of a random $d$-tuple in $\frac{1}{\sqrt{n}} \comhdn$ on an 
open subset of the unit ball in $\R^d$ tend to $0$?
\end{question}
\end{enumerate}

\section{Main results}

In this note, we shall consider random commuting $d$-tuples. 
There are two natural cases to consider, the Hermitian and the general case, which we denote by
$\comhdn$ and $\comdn$.
Since commuting matrices have common eigenvectors, 
the eigenvalues of an element of $\comdn$ are points in $\C^d$, and for an element of $\comhdn$ they are points in $\R^d$. We would like to understand how they are distributed, both for fixed $n$ and as $n$ tends to infinity.

For an alternative approach, where the matrices are not required to commute exactly but are penalized for having large commutators, see the recent paper by Guionnet and Maurel-Segala \cite{gms22}.

\subsection{Main results: Hermitian case}
\label{subsecab}

We only consider weights which are unitarily invariant, i.e. $w(UXU^*) = w(X)$.
We shall call unitarily invariant functions {\em class functions}.
By Proposition \ref{prc1}, with probability $1$ an $X$ in $\comhdn$ will have exactly $n$ distinct eigenvalues.
Since the unitary orbits of $X$ in $\comhdn$ are uniquely determined by their eigenvalues, we can define a 
function $\tw$ on $(\R^d)^n$ by $\tw(\l) = w(X)$, where $X$ is any element of  $\comhdn$ that has $\l$ 
as its eigenvalues.  Writing $\l = \l_j^r$ for $1 \leq j \leq n$ and $1 \le r \le d$, the order of the $j$'s does not matter,
but the order of the $r$'s does.

Our first theorem is a $d$-dimensional analogue of the Ginibre formula.
We write $|\l_i - \l_j |$ to denote the Euclidean distance between $\l_i $ and $\l_j$ in $\R^d$.

{\bf Theorem \ref{thmc2}:} {\em
Let $w:\comhdn \to [0,\i)$  be a continuous integrable class function,
 and assume $\int_{\comhdn} w(X) d X = 1$.
Let random matrices in $\comhdn$ be chosen with distribution $w(X) dX$. 
Then the density of their eigenvalues is given by
\be
\rho_n(\l) \= \ccn \ \tw(\l) \prod_{1 \leq i < j \leq n} |\l_i - \l_j |^2 .
\label{eqc4}
\ee
}

 Note that the repulsion between the eigenvalues is the square of the Euclidean distance in $\R^d$, not of the components.

Using Johansson's large deviation estimate \cite{jo98}, we show that as $n \to \i$ the distribution tends to the equilibrium measure for the logarithmic energy. These have been described explicitly for $d \geq 3$ by
Chafai, Saff and Womersley \cite{csw22, csw23}, and for $d=1$ and $2$ in \cite{wi58} and \cite{ef05}. This yields the following result (see Theorem \ref{thmc5} for a precise statement).
\bt
\label{thmb1}
Let $X_n$ be a random variable in $\comhdn$, and let the weights be probabilities given by
$w_n(X) \= c_n e^{-\frac 12{\rm tr}[ \sum_{r=1}^d (X^r)^* X^r]}$.
As $n \to \i$, the eigenvalues of $\frac{1}{\sqrt{n}}
X_n$ tend to the following equilibrium distribution on the closed ball of radius $R_d$, where $\sigma^d$ is normalized area measure on the sphere, and $d\varrho$ means radial measure.
\begin{eqnarray*}
 \frac{1}{2 \pi } \sqrt{4 - x^2}\ dx, & R_1 = 2& d = 1;
 \\
 \frac{1}{2 \pi } 1_{|x| < \sqrt{2}}\  dx^1 dx^2,
 & R_2 = \sqrt{2} & d =2;
 \\
 \frac{3}{4 \pi^2 } \frac{1}{\sqrt{\frac{4}{3} - |x|^2}} 1_{|x| < R_3}\ d\varrho d\sigma^2,
 &R_3 = \sqrt{\frac{4}{3}} & d = 3;
 \\ \sigma^{d-1}_{R_d},
 & R_d := 1,&d \geq 4.
 \end{eqnarray*}
\et
A corollary of this result (see Corollary \ref{cord1}) is that if one chooses a random commuting Hermitian $d$-tuple in $\comhdn$ and then chooses randomly one of the $d$ matrices in the $d$-tuple, then the normalized eigenvalues of this matrix will tend, as $n \to \i$, to a semi-circular distribution on $[-R_d, R_d]$ for $d \leq 4$ but will be more heavily concentrated near $0$ for any $d \geq 5$.

\subsection{Main results: General case}
\label{subsecac}

In the non-Hermitian case, the set $\comdn$ is irreducible if $d \leq 2$ or $n \leq 3$ or for
a certain range of $n$ when $d=3$ (see Subsection \ref{ssecd1} for details). We restrict our attention to the irreducible case, which happens when the elements with $n$ eigenvalues are generic.
In this case we can write a generic $X \in \comdn$ as
\[
X \= U Q U^* \= U A D A^{-1} U^* ,
\]
where $U$ is a unitary in $\un$, $Q$ is an upper triangular commuting $d$-tuple,
$D$ is a diagonal $d$-tuple and $A$ is an upper triangular unipotent matrix.

In Proposition \ref{prd17} we describe how to take advantage of the unitary invariance of $w_n$ to decompose
$w_n(X) dX$ as $w_n(Q) \kappa(Q) dQ d\nu$, where $\nu$ is the invariant measure on $\un / \T^n$ 
and $\kappa$ is a function described in Subsection \ref{subsecgin}.

In Theorem \ref{thmd2} we describe how to get the eigenvalue distribution of a random element in $\comdn$ by integrating with respect to both $U$ and $A$. 
The formula is unfortunately complicated. In Theorem \ref{thmd1} we calculate the distribution explicitly for commuting $d$-tuples of $2$-by-$2$ matrices. We get:

{\bf Theorem \ref{thmd1}:}
{\em
Let $X$ be a random element of $\vdt$ chosen with distribution
\[
w(X) \= C e^{-\gamma \|X\|_F^2} .
\]
Then the eigenvalues $\l$  of $X$ in $(\C^d)^2$ have distribution
\[
\rho^d(\l) \= 
C'_d \ e^{-\gamma ( |\l_1|^2 + |\l_2|^2)} |\l_2 - \l_1|^2 \sum_{j=0}^{d-1} \frac{1}{(d-1-j)!} \frac{2^j}{\gamma^{j+1}} \frac{1}{|\l_2 - \l_1|^{2j}} .
\]
For $d=1,2,3$ we get respectively
\beq
\rho^1(\l) &\=& C_1^\prime e^{-\gamma ( |\l_1|^2 + |\l_2|^2)}\left[   \frac{1}{\gamma} |\l_2 - \l_1|^2\right]  \\
\rho^2(\l) &\=& C_2^\prime e^{-\gamma ( |\l_1|^2 + |\l_2|^2)} 
\left[ \frac{1}{\gamma} |\l_2 - \l_1|^2 + \frac{2}{\gamma^2}\right] \\
\rho^3(\l)  &\=& C_3^\prime e^{-\gamma ( |\l_1|^2 + |\l_2|^2)} 
\left[ \frac{1}{2\gamma} |\l_2 - \l_1|^2 + \frac{2}{\gamma^2}+ \frac{4}{\gamma^3} |\l_2 - \l_1|^{-2} \right] .
\eeq
}

Notice that, unlike in the Hermitian case or the case $d=1$, the negative powers that appear for $d \geq 3$ show that there is an {\em attraction} between the eigenvalues, not a repulsion.

\section{Known results}
\label{secb}

\subsection{The Gaussian Unitary Ensemble}
\label{secgue}

The Gaussian Unitary Ensemble (or GUE) is the matrix ensemble in which the strictly upper triangular entries
$(\xi_{ij})_{i< j}$ are chosen i.i.d. $N_\C(0,1)$ (ie with a complex normal distribution with mean $0$ and variance $1$), the diagonal entries $\xi_{ii}$ are chosen independently from $N_\R(0,1)$, and the
lower triangular entries are defined by $\xi_{ji} = \overline{\xi_{ij}}$.
This gives an ensemble of Hermitian matrices, parametrized by $\R^{n^2}$. The probability measure
is given by
\be
\label{eqb10}
 C_n e^{-\frac 12  \sum_{i,j=1}^n |X_{ij}|^2 } dm(X) \ee
where $m$ is Lebesgue measure on $\R^{n^2}$, and $C_n$ is the appropriate normalizing constant.

The Ginibre formula \cite{gi65} asserts that the density function of the eigenvalues $\l = (\l_1, \dots , \l_n)$ in $\R^n$ is given by
\be
\label{eqg2}
\rho_n(\l) \= C_n^\prime e^{-|\l|^2/2} \prod_{1 \leq 1 < j \leq n} |\l_i - \l_j |^2 .
\ee
The Wigner semi-circular law  says that if one defines a measure by putting a point mass of weight $1/n$ at each eigenvalue of
$\frac{1}{\sqrt{n}} X$ then these measures  converge weakly to the semi-circular distribution on $[-2,2]$ given by
\[
\frac{1}{2\pi} \sqrt{4-x^2} dx .
\]
See e.g. 
 \cite[Thm. 6.96 and Prop. 6.156]{dei99} 
for a proof.

If instead of choosing each matrix entry independently one chooses Hermitian
matrices with a distribution that replaces \eqref{eqb10} by
\[
C_n e^{-{\rm tr} ( Q(X)) } dX ,
\]
where $Q(x) = \gamma |x|^\alpha + O(|x|^{\alpha -1})$ for some positive constants $\alpha, \gamma$, then the Ginibre formula still holds, becoming
\[
\rho_n(\l) \= C_n^\prime e^{Q(\l)} \prod_{1 \leq 1 < j \leq n} |\l_i - \l_j |^2 
\]
\cite[Thm. 5.22]{dei99}. The semi-circle distribution is replaced by the equilibrium measure for $Q$ \cite[Thm. 6.96]{dei99} (see Subsection \ref{subeq} for a definition of equilibrium measure).

%

\subsection{Equilibrium measures for the logarithmic potential}
\label{subeq}

Let $Q(x)$ be an external field on $\R^d$, a non-negative continous function that goes to infinity 
like some positive power of $|x|$.
For $\mu$ a probability measure on $\R^d$, define the logarithmic energy by
\be
\label{eqbb1}
I^Q(\mu) \= \int \int \log \frac{1}{|x-y|} d\mu (x) d\mu(y) + \int Q(x) d \mu (x) .
\ee
The equilibrium measure for $Q$ is the unique probability measure $\mu^Q$ that minimizes
$I^Q(\mu)$ (for existence and uniqueness, see \cite[Thm. 4.4.14]{bhs19}).

Let us restrict to the case $Q(x) = \gamma |x|^\alpha$.
When $\alpha =2$, the equilibrium measures are known for all $d$.
Let $\sigma^d$ be area measure on the unit sphere $S^d$, normalized to have total mass $1$,
and $\sigma^d_R$ be normalized area measure on the sphere of radius $R$.

\bt
\label{emg}
Let $Q(x) = \gamma |x|^2$. Then the equilibrium measure is given by
\begin{eqnarray}
\label{eqb1}
 \frac{2}{\pi R_1^2} \sqrt{(R_1^2 - x^2)_+}\ dx, &\quad R_1 := \sqrt{\frac 2\gamma}, & d = 1;
 \\
 \frac{1}{\pi R_2^2} 1_{|x| < R_2}\  dx^1 dx^2,
 \label{eqb2}
 &R_2 := \frac{1}{ \sqrt{\gamma}}, & d =2;
 \\
 \frac{1}{\pi^2 R_3^2} \frac{1}{\sqrt{R_3^2 - |x|^2}} 1_{|x| < R_3}\ dr d\sigma^2,
 &R_3 := \sqrt{\frac{2}{3\gamma}},& d = 3;
 \label{eqb3}
 \\ \sigma^{d-1}_{R_d},
 & R_d := \frac{1}{\sqrt{2\gamma}},&d \geq 4.
 \label{eqb4}
 \end{eqnarray}
\et
Result \eqref{eqb1} is the semi-circle law, first obtained by Wigner \cite{wi58}.
For other values of $\alpha$ when $d=1$, see \cite[Sec. 6.7]{dei99}.
Result \eqref{eqb2} is the circular law for random normal matrices; it was proved by Elbau and Felder
in \cite{ef05};
for another proof see \cite{ta12}.
Results \eqref{eqb3} and \eqref{eqb4} are due to Chafai, Saff and Womersley \cite{csw22, csw23}. 
Among other results they also show that the equilibrium for all $\alpha \geq 2$ and all $d \geq 4$ is given by $\sigma^{d-1}_{(1/\gamma \alpha)^{1/\alpha}}$.

\subsection{Other results}

The Hoffman-Weilandt inequality \cite{howi53}
 says that if we let $\mu_j(A)$ denote the $j^{\rm th}$ largest eigenvalue of the self-adjoint matrix $A$, then
 for any self-adjoint $n$-by-$n$ matrices $A$ and $B$ we have
\be
\label{eqab1}
\sum_j | \mu_j (A) - \mu_j(B) |^2 \leq \| A - B \|^2_F .
\ee

\section{Random commuting Hermitian matrices}
\label{secc}

\subsection{Generic $d$-tuples}
\label{subsecc1}

In this section, we shall study random $d$-tuples in $\vhdn$.  
All matrices in this section are Hermitian, and dimension always means real dimension.
First we shall prove that the set of $X$ for which any of the matrices has a repeated eigenvalue is of lower dimension than $\vhdn$, so will have measure $0$.

\begin{definition}
Let $Y$ be a self-adjoint matrix. We shall say $Y$ has banner $(r_1, r_2, \dots , r_p)$, where
$r_1 \geq r_2 \geq \dots \geq r_p \geq 1$ are integers, if $Y$ has $p$ distinct eigenvalues, with multiplcities $r_1, \dots , r_p$.
\end{definition}

Let $\mathcal G(d,n)$ denote the complex Grassmanian of all $d$-dimensional linear subspaces of $\C^n$. It has dimension $2d(n-d)$.
\bprop
\label{prc1}
For $d \geq 2$, the set $\{ X \in \vhdn: \ {\rm banner}(X^1) = (r_1, \dots, r_p) \}$ 
has dimension
  $n^2 + (d-2)n + p$.
\eprop
\bp
We shall prove this by induction on $d$.
First, let $d=1$.
To choose a matrix with banner $(r_1, \dots, r_p)$, one can first choose an $r_1$ dimensional subspace, and the first eigenvalue, then an $r_2$ dimensional subspace of $\C^{n-r_1}$
and the second eigenvalue, and so on. This yields a space of dimension
\begin{multline}
\label{eqc1}
2r_1(n-r_1) + 2 r_2 (n-r_1 - r_2) + \dots + 2 r_{p-1}(n - r_1 - \dots - r_{p-1}) + p\\
\begin{aligned}
 &=\ 2n(n-r_p) - 2 \sum_{i=1}^{p-1} \sum_{j=1}^i r_i r_j + p\\
 &=\ 2n(n-r_p) - \left( \sum_{i=1}^{p-1} r_i \right)^2 - \sum_{i=1}^{p-1} r_i^2 + p\\
 &=\ 2n(n-r_p) - (n-r_p)^2 - \sum_{i=1}^{p} r_i^2 +r_p^2 +  p\\
 &=\ n^2 - \sum_{i=1}^{p} r_i^2 + p.
\end{aligned}
\end{multline}
Note that $n^2 - \sum_{i=1}^{p} r_i^2 + p$ is always $\le n^2$, and equals $n^2$ if each $r_i$ is $1$ (which means $p=n$).

Now assume $d \geq 2$, 
and that we have shown that if $d' < d$, the dimension of
$\vhdpn$ is $n^2 + (d'-1)n$. Note that \eqref{eqc1} proves this assertion for $d' = 1$,
and if $d' > 1$ it follows from the inductive hypothesis, 
since
the maximum value that 
$n^2 + (d'-2)n + p$ can attain is $n^2 + (d'-1)n$, which occurs when $p = n$.
To calculate the dimension of $\{ X \in \vhdn: \ {\rm banner}(X^1) = (r_1, \dots, r_p) \}$, note that once $X^1$ has been chosen, then the dimension of  $\{ X^2,\dots, X^d\}$ 
is $\sum_{i=1}^p {\rm dim}(\vhdmrn)$. By the inductive hypothesis, this gives
\beq
{\rm dim} \{ X \in \vhdn: \ {\rm banner}(X^1) = (r_1, \dots, r_p) \} &\=&
n^2 - \sum_{i=1}^{p} r_i^2 + p + \sum_{i=1}^p( r_i^2 +(d-2) r_i)\\
&=& n^2 +(d-2) n + p.
\eeq
This proves the proposition for $d$, and in particular shows that the dimension of $\vhdn$
is $n^2 + (d-1)n$.
\ep

It follows from Proposition~\ref{prc1} that the set of $X \in \vhdn$ 
for which each $X^r$ has $n$ distinct eigenvalues is of full measure. Let us call these the generic elements, and write
\be
\label{eqc2}
\notag
\comhdng \= \{ X \in \vhdn: \ {\rm each\ }X^r\ {\rm has\ } n\ {\rm distinct\ eigenvalues} \}.
\ee
Note that at any point where the eigenvalues of $X$ are distinct (in particular
at the generic elements) the algebraic set $\vhdn$ is smooth,
 since coordinates can be given by specifying an orthonormal frame for the eigenvectors and a choice of corresponding eigenvalues.

We could choose to order the eigenvalues of $X$, which are points in $\R^n$, by, for example, putting the eigenvalues in dictionary order. However, we will not do this, and 
so the eigenvalues of elements of $\comhdng$ will be invariant under the symmetric group
$\frak{S}(n)$, and when we list them as $(\l_1, \dots, \l_n)$ this is with some random ordering.

\subsection{Ginibre formula}
\label{subsecc2}

Next we prove the $d$-dimensional version of the Ginibre formula \eqref{eqg2}.
Both the proof in \cite[Sec. 2.6.1]{ta12} and in
\cite[Thm. 5.22]{dei99} of the one dimensional case can be adapted to $d$ dimensions, though some care must be
 taken since the tangent space to $\vhdn$ can no longer be naturally identified with $\vhdn$.
In the $d$-variable case, the eigenvalues of an elements of $\vhdn$ will be points in
$(\R^d)^n$, and we will write them as $\l = (\l_1, \dots, \l_n)$ where each $\l_i$ is a $d$-tuple
$\l_i = (\l_i^1, \dots, \l_i^d)$.

If $w$ is a  real-valued function on $\comhdng$ that is invariant under the action of the unitary group $\un$  (ie. $w(X) = w(U^* X U)$ for all $U \in \un$) then $w(X)$ 
can only depend on the eigenvalues of $X$, since they completely determine the unitary equivalence class of $X$. Such functions are called {\em class functions.}
If $w$ is a class function, there
is a  function $ \tw : (\R^d)^n \to \R$ that is invariant with respect to the action of
 $\frak{S}(n)$ so that
 \be
 \label{eqc3}
 \notag
 w(X) \= \tw (\l),
 \ee
 where $\l =  (\l_1, \dots, \l_n)$ is some ordering of the eigenvalues of $X$.
 Moreover, if $w$ is continuous, so is $\tw$.

For $x \in \R^d$, we shall use $|x|$ to mean the Euclidean norm $\sqrt{\sum_{r=1}^d |x^r|^2}$.
We shall write $dX$ to mean $d {\mathcal H}_{{\rm dim} \vhdn} (X)$.
We shall use the subscript $F$ for the Frobenius (or Hilbert-Schmidt) norm on $\vdn$, viz.
\[
\| X \|_F^2 \= \sum_{r=1}^d \sum_{i,j =1}^n |X_{ij}^r|^2 .
\]
Let $v(k)$ denote the volume of the unit ball in $\R^k$.

Conjugating a diagonal $d$-tuple by a unitary gives the most general Hermitian $d$-tuple, but
diagonal unitaries have no effect, so we must understand $\un / \T^n$, where we identify a point $t \in \T^n$,
the set of $n$-tuples of unimodular complex numbers, with the diagonal unitary $D_t$ that has the entries of $t$ on the diagonal.
 We shall write $1_n$ for the identity
in $\T^n$.
The tangent space at the identity of $\un$ is the Lie algebra $\frak{u}(n)$ which consists of the skew-Hermitian matrices (matrices satisfying $S^* = - S$). Let $\azn$ denote the subspace of skew-Hermitian matrices that are $0$ on the diagonal.

\begin{lemma}
\label{lemc11}
Let $g: \azn \times \T^n \to \un$ be defined by $g(S,t) = e^S D_t$. Then $dg : \azn \times \R^n
\to \frak{u}(n)$ is an isomorphism between the tangent space of $g: \azn \times \T^n$
at $(0,1_n)$ and the tangent space of $\un$ at the identity.
\end{lemma}
\bp
Note that
\[
dg|_{(0,1_n)} : ( S, \theta) \ \mapsto \ S + D_{i\theta} ,
\]
and that this map is an isomorphism.
\ep
In particular, by the inverse function theorem one can parametrize a neighborhood of $[\T^n]$ in
$\un/\T^n$ as $\{ e^S \}$ where $S$ ranges over a neighborhood of $0$ in $\azn$.
Define a constant $\ccn$ by
\be
\label{eqc21}
\ccn \= \lim_{\vare \to 0^+} 
\frac{ {\rm vol}_{\azn}\{ S :  S \in \azn, \| S \|_F < \vare \}}
{ {\rm vol}_{\un/\T^n} \{ e^S : S \in \azn, \| S \|_F < \vare\} }
 .
\ee
The number $\frac{1}{\ccn}$  is the Jacobian determinant of $e^S$, normalized by the requirement that the total Haar measure of $\un$ is $1$.
\begin{lemma}
\label{lemc12}
Let $X \in \comhdng$. The tangent space $T_X \vhdn$ is the set $Z = (Z^1, \dots, Z^d)$ 
of  $d$-tuples of self-adjoint matrices that satisfy 
\be
\label{eqc233}
[X^r , Z^s ] \= [X^s, Z^r] \qquad \forall \ 1 \leq r \neq s \leq d .
\ee
If $e_j$ is an orthonormal basis of eigenvectors of $X$ with corresponding eigenvalues $\l_j$,
then with respect to this basis equations \eqref{eqc233} become
\be
\label{eqc24}
Z^s_{ij} (\l_j^r - \l_i^r) \= Z^r_{ij} (\l_j^s - \l_i^s)\qquad \forall \ 1 \leq r \neq s \leq d .
\ee
\end{lemma}
This result has been known for a long time---we do not know where it first appeared.
For the convenience of the reader, we include the proof from 
\cite[Lemma 7.3]{amy12c}. (One could also
 adapt the proof of Lemma \ref{lemd1} below).

\bp
We want to show that  there exists a $C^1$ curve $X(t)$ of commuting self-adjoints
with $X(0) = X$ and
$X'(0) = Z$ if and only if  \eqref{eqc24} holds.

($\Rightarrow$): If $X(t) = X + t Z + o(t)$ is commutative, calculate
$$
[X^r (t) , X^s(t)] \= t \left( [ X^r, Z^s] - [ X^s , Z^r] \right) + o(t).
$$
The coefficient of $t$  must vanish, giving (\ref{eqc24}).

($\Leftarrow$):
Suppose 
(\ref{eqc24}) holds. This means
\be
\notag
Z^s_{ij} ( \lambda^r_j - \lambda^r_i ) \=
Z^r_{ij} ( \lambda^s_j - \lambda^s_i ) \qquad \forall \, r \neq s ,
\ee
so
\be
\label{eqm06}
Z^r_{ij} \frac{1}{\lambda^r_j - \lambda^r_i} \=
Z^s_{ij} \frac{1}{\lambda^s_j - \lambda^s_i}  \qquad {\rm if\ }
\lambda^r_j - \lambda^r_i \neq 0 \neq
\lambda^s_j - \lambda^s_i.
\ee
Define a skew-selfadjoint matrix $Y$ by
\be
\label{eqm07}
Y_{ij} \= 
Z^r_{ij} \frac{1}{\lambda^r_j - \lambda^r_i} \qquad {\rm for\ any\ }r\ {\rm such\ that\ }
\lambda^r_j - \lambda^r_i \neq 0 .
\ee
For any $i \neq j$, there is some $r$ with 
$
\lambda^r_j - \lambda^r_i \neq 0$, so (\ref{eqm07}) defines $Y_{ij}$;
and (\ref{eqm06})
says it doesn't matter which $r$ we choose.
Let all the diagonal terms of $Y$ be $0$.

Define
\be
\notag
X^r(t) \= e^{tY}( X^r +  t\ {\rm diag}( Z^r))  e^{-tY} .
\ee
Since $e^{tY}$ is a unitary matrix and $X^r  + t\ {\rm diag}( Z^r)$ is diagonal, $X(t) \in CXAM_n^d$ and
\[
 \frac{d}{dt}X^r(t)|_{t=0} = [Y,X^r] + {\rm diag}( Z^r )= Z^r.
\]
\ep

\begin{theorem}
\label{thmc2}
Let $w:\vhdn \to [0,\i)$  be a continuous integrable class function,
 and assume $\int_{\vhdn} w(X) d X = 1$.
Let random matrices in $\vhdn$ be chosen with distribution $w(X) dX$. 
Then the density of their eigenvalues is given by
\be
\rho_n(\l) \= \ccn \ \tw(\l) \prod_{1 \leq i < j \leq n} |\l_i - \l_j |^2 .
\label{eqc4x}
\ee
where $\ccn$ is given by \eqref{eqc21}.
\et
The theorem asserts that if we wish to integrate some class function  with respect
to $w(X) dX$, we can first integrate over $\un/\T^n$, and then integrate over the eigenvalues.
This gives the decomposition
\be
\label{eqc45}
w(X) d X \= \rho_n(\l) d\nu d\l 
\ee
where $d\l$ is Lebesgue measure on $(\R^d)^n$ and $d\nu$ is volume measure on $\un/\T^n$, with
$\rho_n$ given by \eqref{eqc4x}. Notice that the repulsion between the eigenvalues is based on their Euclidean separation in $\R^d$.

Let $E_{ij}$ denote the matrix which is $1$ in the $(i,j)$ entry and $0$ elsewhere.
\bp
The idea is to decompose $w(X) dX$ as in \eqref{eqc45}. 

As every element in $\vhdn$ is unitarily equivalent to a diagonal tuple, it is sufficient to look at
a neighborhood of a diagonal tuple.
Fix some orthonormal basis of $\C^n$ and choose a diagonal tuple $D$ in $\comhdng$ with diagonal entries $\l = (\l_1, \dots , \l_n)$ in $(\R^d)^n$. Let $\vare$ be a sufficiently small positive number,
and consider $\{ X : \| X - D \|_F < \vare \}$.

By the  Hoffman-Wielandt inequality  \eqref{eqab1}, the eigenvalues of $X$ 
can be ordered as $\kappa = (\kappa_1, \dots, \kappa_n)$
where $\kappa_j$ is the eigenvalue of $X$ closest to $\l_j$, and
\be
\sum_{r=1}^d \sum_{j=1}^n | \kappa_j^r - \l_j^r |^2 \ \leq \ \| X - D \|^2_F  \ < \ \vare^2 .
\label{eqc65}
\notag
\ee
By Lemma \ref{lemc11}, we can parametrize a small neighborhood of $D$ by
\[
X \= G(\theta, S) \= e^S ( D + D_\theta) e^{-S}
\]
where $S$ ranges over a neighborhood of $0$ in $\azn$ and $\theta$ over a neighborhood
of $0$ in $\rdn$. The derivative $dG$ at $(0,0)$  is the map between the tangent spaces
\beq
dG : \rdn \times \azn   & \ \to \ & T_{D} \vhdn \\
 (\theta, S) &\ \mapsto \ &SD - DS + D_\theta .
\eeq
To calculate the Jacobian determinant of $G$,
let $R_{ij} = \frac{1}{\sqrt{2}} (E_{ij} - E_{ji})$ and
$S_{ij} = \frac{i}{\sqrt{2}} (E_{ij} + E_{ji})$.
Then we can use the real coordinates $\theta$ for $\rdn$ 
and $\{ R_{i,j}, S_{ij} : i < j \}$ for $\azn$.
We have
\beq
dG: (\theta,0) &\ \mapsto \ & \theta\\
dG: (0, R_{ij}) &\mapsto& -i  (\l_j^r - \l_i^r) S_{ij} \\
dG :(0, S_{ij}) & \mapsto & i(\l_j^r - \l_i^r) R_{ij} .
\eeq
Therefore the Jacobian of $G$ is
\[
\sqrt{ |\det (dG^* dG) |} \=
 \prod_{1 \leq i < j \leq n} |\l_i - \l_j |^2  .
 \]
 This gives 
 \[
 dX \=
  \prod_{1 \leq i < j \leq n} |\l_i - \l_j |^2 \ d\theta d S.
 \]
By Lemma \ref{lemc11}, $dS$ is locally $\ccn d\nu$.
Thus we have shown that 
\[
w(X) dX =  \ccn \tw(\l)  \prod_{1 \leq i < j \leq n} |\l_i - \l_j |^2  d\l d \nu ,
\]
which gives \eqref{eqc4x}.
\ep

\begin{remark}{\rm
As the constant $\ccn$ depends only on the unitary group $\un$ and not on $d$ or $w$,
it can be calculated in a number of ways; see e.g.  \cite{gi65, ta12}.
}\end{remark}

\subsection{Density of eigenvalues}

In Theorem \ref{thmc2} we gave a formula for the joint distribution of eigenvalues in $(\R^d)^n$. Now we ask what is the density of eigenvalues in $\R^d$, as $n \to \infty$. This will depend on the equilibrium measure of $ \log \tw$, and choosing the correct scaling.
We shall restrict ourselves to weights of the following form, where $\alpha, \gamma > 0$:
\be
\label{eqc10}
Q(x) \= \gamma (\sum_{r=1}^d |x^r|^2)^{\alpha/2}.
\ee
The following theorem is a special case of \cite[Thm. 4.4.14]{bhs19}.
\bt
\label{thmc3}
Let $Q$ be as in \eqref{eqc10}, and let the logarithmic energy be defined as in
\eqref{eqbb1}. Then there is a unique probability measure $\mu^Q$ on $\R^d$ where $\inf I^Q(\mu)$ 
is attained, and moreover $\mu^Q$ is compactly supported.
\et
We call $\mu^Q$ the equilibrium measure (for $Q$) and set
\[
E^Q \= I^Q(\mu^Q) .
\]

We define $w_n$ on $\vhdn$ by
\be
w_n(X) \= C_n \ e^{-{\rm tr}\ Q(X)},
\label{eqc11}
\notag
\ee
where $C_n$ is a normalizing constant to ensure 
$\int_{\vhdn} w_n(X) dX = 1$.
If $X \in \comhdng$ and $\l$ is some ordering of the eigenvalues, then
\be
\notag
\tw_n(\l) \= C_n e^{- \sum_{j=1}^n Q(\l_j)} .
\ee

Let $\rho = \rho_n$ be as in \eqref{eqc4x}, so
\be
\label{eqc12}
\notag
\rho_n(x) \= \ccn C_n \, \,  e^{- \sum_{j=1}^n Q(x_j)} \prod_{1 \leq i < j \leq n} |x_i - x_j|^2 .
\ee
 To find the distribution of a single eigenvalue, we must integrate
out the other $n-1$ eigenvalues. First we rescale by $ y = n^{-\frac{1}{\alpha}} x$, defining a new probability density on 
$(\R^d)^n$ by
\be
\label{eqc13}
\cp_n(y) \= D_n\,  e^{- n \sum_{j=1}^n Q(y_j)} \prod_{1 \leq i < j \leq n} |y_i - y_j|^2 .
\ee
Here $D_n$ is a normalizing constant to make $\int_{(\R^d)^n} \cp_n(y) dy = 1$.
It is related to the previous constants by
\[
D_n \= \left(n^{\frac{1}{\alpha}}\right)^{n^2 + (d-1)n} \ \ccn C_n .
\]
If $\phi$ is a bounded continuous function on $\R^d$, we define an expectation by
\[
\en( \phi) \ := \ \int_{\rdn} \frac 1n \sum_{i=1}^n \phi(y_i) \cp_n(y) dy .
\]
The principal result of this section is the following theorem.
\bt
\label{thmc6}
Let $Q$ be as in \ref{eqc10}, and let $\phi$ be a bounded and continuous function on $\R^d$. Then
\be
\label{eqc18}
\lim_{n \to \i} \en(\phi) \=
\int_{\rd} \phi(x) d\mu^Q(x) .
\ee
\et

Theorem \ref{thmc6} was proved 
 in the $d=1$ case by K. Johansson \cite{jo98},
 and his proof generalizes to the $d$-dimensional case.
 For the convenience of the reader we shall give a complete proof, following the exposition 
 in Section 6.4 of the beautiful monograph \cite{dei99}, but we emphasize that generalizing to $d$ dimensions
 requires little more than notational changes.

We start with a large deviation estimate, which for $d=1$ is Lemma 4.2 in \cite{jo98}
or Lemma 6.67 in \cite{dei99}.
Our proof is a straightforward generalization, but since the proof that Claim \ref{eqc17} holds
must be modified for $d \geq 2$,  we include the whole proof for completeness.

Define, for $s,t \in \rd$, 
\[
k(s,t) \ := \ \log \frac{1}{|s-t|} + \frac 12 ( Q(s) + Q(t) ),
\]
and 
 for $y = (y_1, \dots, y_n)  \in \rdn$,
\begin{eqnarray}
\label{eqc15}
\notag
K_n(y) &\= & \sum_{i \neq j} k(y_i, y_j) \\
\notag
\\&=& \sum_{1 \leq i < j \leq n} \log |y_i - y_j|^{-2} + (n-1) \sum_{i=1}^n Q(y_i) .
\end{eqnarray}
Let
\be
\notag
A_{n,\eta} \ :=\ \{ y \in (\R^d)^n : \frac{1}{n^2} K_n(y) \leq E^Q + \eta \} .
\ee

\begin{lemma}
\label{lemc1}
Let $\cp_n(y)$ be the probability density given by \eqref{eqc13}.  
Let $\eta, a > 0$.
Then there exists $n_0$ which depends on $\eta$ but not on $a$ so that 
\be
\label{eqc16}
\notag
\int_{{(\R^d)^n} \setminus A_{n, \eta + a}} 1 \ d\cp_n  \
 \leq\  e^{-an^2}, \quad n \geq n_0.
\ee
\end{lemma}
\bp
Assume $d \geq 2$, as the $d=1$ case is proved in \cite{jo98}.
Let $ 0 < \vare < \eta/2$.
For $\delta > 0$, let
\[
\psi_{\delta} (y) \= \frac{1}{v(d)\delta^d} \int_{B(y,\delta)} d\mu^Q(x) ,
\]
where $v(d)$ denotes the volume of the unit ball in $\R^d$, and $B(y,\delta)$ means the ball centered at $y$ of radius $\delta$.
Then $\psi_{\delta}$ is a compactly supported function, and 
$\psi_{\delta}(y) dy$ tends vaguely to $\mu^Q$ as $\delta \to 0$.
Warning: $\psi_\delta$ need not be continuous, for example if $\mu^Q$ puts mass on a sphere of radius $\delta$.

 Claim: For $\delta$ sufficiently small, we have
\be
\label{eqc17}
I^Q (\psi_{\delta}) \ \leq \ E^Q + \frac{\vare}{2}.
\ee
{\sc Proof of Claim:} 
\beq
\lefteqn{\int_{\R^d} \int _{\R^d} \log |x-y| \psi_{\delta}(x) \psi_{\delta}(y) dx dy }\\
&\=&
\int \int \left[ \frac{1}{v(d)^2\delta^{2d}} \int_{B(0,\delta)}   \int_{B(0,\delta)} 
\log |x-y + s-t | \ ds dt
  \right] d\mu^Q(x) d\mu^Q(y).
  \eeq
  Since $d \geq 2$, the function $\log |x|$ is subharmonic. 
  Applying subharmonicity of $\log|x|$ twice, we get
\[
  \frac{1}{v(d)^2\delta^{2d}} \int_{B(0,\delta)}   \int_{B(0,\delta)} 
\log |x-y + s-t |  \ ds dt
\ \geq \ \log|x-y|.
\]
Therefore we conclude
\[
\int_{\R^d} \int _{\R^d} \log \frac{1}{|x-y|} \psi_{\delta}(x) \psi_{\delta}(y) dx dy 
\ \leq \ 
\int_{\R^d} \int _{\R^d} \log \frac{1}{|x-y|} d\mu^Q(x) d\mu^Q(y).
\]
As
\[
\lim_{\delta \to 0} \int Q(x) \psi_{\delta}(x) dx \= \int Q(x) d\mu^Q(x) ,
\]
we have proved the claim.

\vs
The rest of the proof proceeds exactly as in \cite{dei99}.
Fix $\delta$ so that \eqref{eqc17} holds. Let $E$ be the compact support of $\psi_{\delta}$. 
We have
\beq
\frac{1}{D_n}
&\=& 
\int_\rdn e^{-K_n(y) - \sum_{i=1}^n Q(y_i)} dy \\
&\geq &
\int_{E^n} e^{-K_n(y) - \sum_{i=1}^n Q(y_i) - \sum_{i=1}^n \log \psi_{\delta}(y_i)}
\prod_{i=1}^n \psi_{\delta}(y_i) dy.
\eeq
Fix some positive $\delta$ for which \eqref{eqc17} holds.
Applying Jensen's inequality to the probability measure $\prod_{i=1}^n \psi_{\delta}(y_i) dy$ we get
\beq
\log D_n &\ \leq \ &
\int \left(  K_n(y)  +  \sum_{i=1}^n Q(y_i) 
+  \sum_{i=1}^n \log \psi_{\delta}(y_i) \right) \ \prod_{i=1}^n \psi_{\delta}(y_i) dy\\
&=& 
n(n-1) I^Q ( \psi_{\delta}) + n \int_{\R^d} Q(x) dx + n \int _{\R^d} (\log \psi_{\delta}(x)) \psi_{\delta}(x) dx \\
&\leq & n^2 I^Q ( \psi_{\delta}) + O(n) \\
&\leq& n^2( E^Q + \frac{\vare}{2}) + O(n).
\eeq
Therefore there exists $n_1$ so that if $n \geq n_1$, then
\[
\frac{1}{n^2} \log D_n \ \leq \ E^Q + \vare .
\]
So for $n \geq n_1$, and recalling that $\vare < \eta/2$, we have
\beq
\int_{ \rdn \setminus A_{n, \eta + a} } 1 \ dy &\=&
D_n \int_{\{ K_n(y) \geq n^2( E^Q + \eta + a\} } e^{-K_n(y) -\sum Q(y_i) } dy \\
&\leq & 
\int_{\rdn} e^{  -\sum Q(y_i)} e^{n^2(E^Q + \vare)} e^{-n^2(E^Q + \eta + a)} dy\\
&\leq&
\left[\left( \int_{\R^d} e^{-Q(x)}dx \right)^n e^{- \eta n^2/2} \right] e^{-an^2} .
\eeq
Now choose $n_0 \geq n_1$ so that the term in brackets is bounded by $1$, and the lemma is proved.
\ep

With Lemma \ref{lemc1} in hand, the proof of Johansson's theorem \cite[Thm. 2.1]{jo98} generalizes to $\R^d$. In particular, the proof of \cite[Thm. 6.96]{dei99} applies with only notational changes.

 {\em Proof  of Theorem \ref{thmc6}.}
 First, we will show that
 \be
 \label{eqc20}
 \limsup_{n \to \i} \en(\phi) \ \le \ 
\int_{\rd} \phi(x) d\mu^Q(x) .
\ee
 Fix $\eta > 0$. Since $\phi$ is bounded, it follows from Lemma \ref{lemc1} that the left-hand side of \eqref{eqc20}
 is the same as
 \[
  \limsup_{n \to \i} \int_{\ant} \frac 1n \sum_{i=1}^n \phi(y_i) \cp_n(y) dy .
  \]
Let $y* = (y_1^*, \dots , y_n^*)$ be a point where the function $\frac 1n \sum_{i=1}^n \phi(y_i)$ attains its
maximum on the compact set $\ant$.
Let $\nu_n$ be the probability measure on $\rd$ that has mass $\frac 1n$ at each point $y_i^*$.
We have
\be
\label{eqc24x}
\notag
\int_{\ant} \frac 1n \sum_{i=1}^n \phi(y_i) \cp_n(y) dy
\ \leq \
\int_{\rd} \phi(t) d \nu_n(t) .
\ee

If $s,t \in \rd$ then $|s-t| \leq \sqrt{(1+|s|^2)(1+|t|^2)}$, so
\beq
(n-1) \sum_{i=1}^n Q(y_i^*) - (n-1) \sum_{i=1}^n \sum \log(1 +|y_i^*|^2) 
&\ \leq \ & K_n(y^*) \\
&\leq& n^2( E^Q + 2\eta) .
\eeq
Therefore
\be
\label{eqc21x}
\int_{\rd} Q(t) - \log (1 +|t|^2) \ d\nu_n(t)
\ \leq \ \frac{n}{n-1} (E^Q + 2 \eta) .
\ee
It foilows from \eqref{eqc21x} that the sequence $\nu_n$ is tight.
(Recall that a sequence of probability measures $\nu_n$ is called tight if for every $\vare > 0$ there exists a compact
set whose complement has measure less than $\vare$ for all sufficiently large $n$).
Since the sequence is tight, there is some subsequence $\nu_{n_j}$ that converges weakly to some probability measure
$\nu^\eta$ for which
\[
\int \phi d \nu^\eta \= \limsup_{n \to \i} \int \phi d \nu_n .
\]
(First extract a subsequence whose limit is the limsup, then extract a further subsequence that converges weakly).

Observe that for any $L > 0$ we have
\beq
\frac{1}{n^2} K_n(y^*) &\=& \frac{1}{n^2} \sum_{i \neq j} k(y_i^*, y_j^*) \\
&\ge& \frac{1}{n^2} \left[ \sum_{i,j=1}^n \min( L, k(y_i^*, y_j^*) - n L \right] \\
&=&
\int \int \min(L, k(s,t)) d\nu_n(s) d\nu_n(t) - \frac{L}{n}  .
\eeq
Since $y^* \in \ant$, we get
\be
\label{eqc22}
\int \int \min(L, k(s,t)) d\nu_n(s) d\nu_n(t) 
\ \leq \ \frac Ln + E^Q + 2 \eta .
\ee
Apply \eqref{eqc22} with $n = n_j$, let $j \to \infty$ and then let $L \to 0$.
This gives
\be
\label{eqc23}
\int \int k(s,t) d\nu^\eta (s) d\nu^\eta(t) \ \leq \ E^Q + 2 \eta .
\ee
Now let $\eta_j$ be a sequence that decreases to $0$.
From \eqref{eqc21x} we see that $\nu^{\eta_j}$ is tight, so it has a subsequence that converges weakly to some measure
$\nu$. 
From \eqref{eqc23}, we see that $\nu$ must be the equilibrium measure $\mu^Q$.
So we have
\be
\label{eqc25}
\notag
 \limsup_{n \to \i} \en(\phi) \ \le \  
 \int \phi \ d\nu^{\eta_j}
 \ee
 for all $j$, and letting $j \to \infty $ we get \eqref{eqc20}.
 
 Repeating the argument with $y_*$ where $\frac 1n \sum \phi(y_i)$ attains its minimum, we get the reverse inequality 
 for the liminf, and hence we get \eqref{eqc18}.

\subsection{Gaussian Case}

 Let $Q(x) = \gamma |x|^2$. 
 (The positive constant $\gamma$ is just a simple scaling, and plays no important r\^ole.
But as some authors default to $\gamma =\frac 12$ and others to $\gamma = 1$, we include it to make the formulas seem more familiar.)
  Let $X_n$ be chosen in $\vhdn$  with distribution 
$C_n e^{-{\rm tr}(Q(X))}$, with $C_n$ the normalizing constant. 
Let $E_n$ denote expectation with respect to this probability.
Then for $\phi$ a continuous bounded function on $\R^d$, we have
\beq
E _n\left[ \frac 1n\  {\rm tr} \, (\phi(\frac{1}{\sqrt{n}} \ X_n)) \right]
&\=&
\int_{\vhdn} \frac 1n\ {\rm tr}\, ( \phi(\frac{1}{\sqrt{n}} X)) w_n(X) dX\\
&=&
\frac 1n \int_{\rdn} \sum_{i=1}^n \phi(\frac{1}{\sqrt{n}} \l_i) \rho_n(\l) d\l \\
&=& 
\frac 1n \int_{\rdn} \sum_{i=1}^n \phi(y_i) \cp_n(y) dy \\
&=& \int_{\R^d} \phi(x) R_{1,n} (x) d x .
\eeq

So combining Theorem \ref{thmc6} with Theorem \ref{emg}, we get the following result in the Gaussian case.
\bt
\label{thmc5} Let $Q(x) = \gamma |x|^2$.  Let $X_n$ be chosen in $\vhdn$  with distribution 
$C_n e^{-{\rm tr}(Q(X))}$. Then for all $\phi$ that are continuous and bounded on $\R^d$, 
we have
\[
\lim_{n \to \i} E_n \left[ \frac 1n\  {\rm tr} \, (\phi(\frac{1}{\sqrt{n}} \ X_n)) \right]
\=
\int_{\R^d} \phi(x) d\mu^Q(x) ,
\]
where the equilibrium measure $\mu^Q$ is given by Theorem \ref{emg}.
\et

Let us consider what happens 
when we  project 
 the equilibrium measures for the Gaussian case 
onto the $x_1$ axis.
%
For $d=2$, the density is
\[
\frac{1}{\pi R_2^2}
\int_{-\sqrt{R_2^2 - x_1^2}}^{\sqrt{R_2^2 - x_1^2}} 1 dx_2 \= \frac{2}{\pi R_2^2} \sqrt{R_2^2 - x_1^2} .
\]
For $d=3$, it is the weighted area of a disk of radius $\sqrt{R_3^2 - x_1^2}$, where the weight
at radius $s$ is $\frac{1}{\pi^2 R_3^2}\frac{1}{\sqrt{R_3^2 - x_1^2 - s^2}}$.
Switching to polar coordinates on this disk, the density is
\[
\frac{2}{\pi R_3^2} \int_0^{\sqrt{R_3^2 - x_1^2}} \frac{s}{\sqrt{R_3^2 - x_1^2 - s^2}} ds \=
\frac{2}{\pi R_3^2} \sqrt{R_3^2 - x_1^2}.
\]
For $d \geq 4$,  we project normalized area measure for the $(d-1)$ sphere $S^{d-1}$ with radius $R_{d}$.
Let us change variables to  $x_1 = R_d \cos \theta$ for $0 \leq \theta \leq \pi$, where $\theta$
represents the angle with the positive $x_1$ axis.
The area of the shell between $\theta$ and $\theta+\Delta \theta$ is approximately
$\Delta \theta$ times the $(d-2)$-dimensional area of a sphere of radius $\sin \theta$,
which is $c R_d^{d-2} \sin^{d-2} \theta$. As $\Delta x = - R_d \sin \theta \Delta \theta$, we get that the 
area of the shell between $x$ and $x + \Delta x$ is approximately $\Delta x$ times
\[
c R_d^{d-3} \sin^{d-3} \theta \= c \left(\sqrt{R_d^2-x_1^2} \right)^{d-3} .
\]
Letting $\Delta x \to 0$, and plugging in the
normalizing value for $c$ which comes from the identity
\[
\int_{-1}^1 (\sqrt{1-x^2})^{d-3} dx \= \frac{\sqrt{\pi}\ \Gamma(\frac{d-1}{2})}{\Gamma(\frac d2)} ,
\]
which is valid for $d > 1$, we get that the density is
\be
\label{eqc21b}
\notag
 \frac{ \Gamma(\frac d2)}{\sqrt{\pi} R_d^{d-2}\ \Gamma(\frac{d-1}{2})} \ \left(\sqrt{R_d^2-x_1^2} \right)^{d-3}.
\ee
For $d = 2,3,4$ we still get a semi-circular distribution, though the radius $R_d$ shrinks as $d$ increases from $1$ to $4$. For $d > 4$ the radius $R_d$ stays $\frac{1}{\sqrt{2\gamma}}$, but the eigenvalues of $X^1$
cluster more around $0$. We have proved the following.
\begin{corollary}
\label{cord1}
 Let $X_n$ be chosen in $\vhdn$  with distribution 
$C_n e^{- \gamma \,{\rm tr}(\sum_{j=1}^d (X^j)^2)}$. Then for all $\phi$ that are continuous and bounded on $\R$, 
we have
\[
\lim_{n \to \i} E_n \left[ \frac 1n\  {\rm tr} \, (\phi(\frac{1}{\sqrt{n}} \ X^1_n)) \right]
\= \int_{-R_d}^{R_d}\phi(x)  f_d(x) dx,
\]
where $R_d$ are as in Theorem \ref{emg}
and
\beq
f_1(x) &\=& \frac{2}{\pi R_1^2} \sqrt{R_1^2- x^2} \\
f_2(x) &=& \frac{2}{\pi R_2^2} \sqrt{R_2^2 -x^2} \\
f_3(x) &\=& \frac{2}{\pi R_3^2} \sqrt{R_3^2 - x^2}\\
f_d(x) &\= & \frac{ \Gamma(\frac d2)}{\sqrt{\pi} R_d^{d-2}\ \Gamma(\frac{d-1}{2})} \ \left(\sqrt{R_d^2-x^2} \right)^{d-3}, \qquad d \geq 4.
\eeq
\end{corollary}
%

\section{Random Commuting Matrices: the non-Hermitian case}
\label{secd}

\subsection{Generic $d$-tuples}
\label{ssecd1}

In this section, we shall drop the assumption that the matrices are Hermitian, and work with random matrices from $\vdn$. 

There is a qualitative difference between the algebraic sets $\vdn$ for different values of $d, n$.
In particular, there are major differences between the reducible and irreducible cases.
Motzkin and Taussky proved that $\vtn$ is irreducible and of dimension $n^2 + n$ for all $n$ \cite{mt55}.
%
Gerstenhaber reproved this, and proved that ${\frak C}^4_4$ is reducible \cite{ge61}.
Guralnick proved that $\vdn$ is irreducible for all $d$ if $n \leq 3$ \cite{gu92}.

When $d=3$, Guralnick proved that ${\frak C}^3_{32}$ is reducible \cite{gu92}.
It is now known that ${\frak C}^3_n$ is reducible if $n \geq 29$ \cite{holom01,ng14},
and irreducible if $n \leq 10$ (see \cite{gs00,holom01,om04,han05, si08,si12b}  for
$n$ going from $4$ to $10$.)

Let us call $X \in \vdn$ {\em generically diagonalizable} if each component $X^r$ has $n$ distinct eigenvalues, and define
\be
\label{eqd1}
\notag
\comdng \= \{ X \in \vdn: \ {\rm each\ }X^r\ {\rm has\ } n\ {\rm distinct\ eigenvalues} \}.
\ee

When are the generically diagonalizable $d$-tuples dense in $\vdn$?
In \cite{holom01}, Holbrook and Omladic showed that the answer is yes 
if and only if $\vdn$ is irreducible. In the irreducible case, the 
non-diagonalizable elements are of smaller dimension than the diagonalizable ones.  Guralnick proved that ${\rm dim} \ {\frak C}^3_n > {\rm dim}\  {\frak C}^3_{n,{\rm gen}}$ if $n \geq 32$.
For a study of the components in the non-irreducible case, see \cite{js22}.
\begin{proposition}
\label{prd1}
Let $\vdn$ be irreducible. Then
\be
\label{eqd2}
{\rm dim}_\C (\vdn \setminus \comdng)  \ < \ 
{\rm dim}_\C (\comdng) = n^2 +(d-1)n .
\ee
\end{proposition}
\bp
Let us write ${\mathcal W}^r = \{ X \in \vdn : X^r {\rm \ does\ not\ have\ }n\ {\rm distinct\ eigenvalues} \}$.
Let $q^r$ be the discriminant of the characteristic polynomial of $X^r$, thought of as a polynomial
in the $dn^2$ entries of an element of $\mn^d$, and let $q = q^1 \cdots q^d$. Then
\[
{\mathcal W}^r \= \{ X \in \vdn: q^r(X) = 0 \} ,
\]
and
\[
\vdn \setminus \comdng \= \cup_{r=1}^d {\mathcal W}^r \= \{X \in \vdn : q(X) = 0 \}.
\]
As $q$ does not vanish identically on $\vdn$, assumed irreducible,  it follows that
\[
{\rm dim}_\C (\vdn \setminus \comdng)  \ < \ 
{\rm dim}_\C (\vdn)  .
\]

When $\vdn$ is irreducible, it is shown in \cite{gs00} that ${\rm dim}_\C(\vdn) = n^2 + (d-1)n$.
It is obvious that ${\rm dim}_\C (\comdng) = n^2 +(d-1)n $, so we get \eqref{eqd2}.
\ep

\subsection{Ginibre formula}
\label{subsecgin}

Throughout the remainder of Section \ref{secd}
 we shall assume that $\vdn$ is irreducible.
\begin{center}
{\em Assumption I: $\vdn$ is irreducible.}
\end{center}
(If this is not the case, then we must replace $\vdn$ by its principal component, that is the closure of the set of jointly diagonalizable matrices).
Let  $w: \vdn \to [0,\i)$ be a continuous function that satisfies $w(U^*XU) = w(X)$ for all unitary $U \in \un$, and is normalized so
\[
\int_{\vdn} w(X) d X \= 1 ,
\]
where $dX$ is Hausdorff measure of the appropriate dimension 
restricted to $\vdn$ (from \eqref{eqd2} the appropriate dimension is $2n^2 + 2(d-1) n$).

Let $X \in \vdn$. There is a unitary $U$ in  $\un$ so that the $d$-tuple $U^* X U$ is jointly upper
triangular. Once such a $U$ has been chosen, let $\l_j \in \C^d$ be the point whose $r^{\rm th}$ component is the $j^{\rm th}$ diagonal entry of $U^* X^r U$.
We shall call the multi-set $\{ \l_1, \dots , \l_n \}$ the {\em multi-spectrum} of $X$.
It is straightforward to see that it is independent of the particular choice of $U$, that it consists of
the eigenvalues of $X$, and that the multiplicity of a particular eigenvalue $\l$ is
\[
\dimc \left( \cap_{r=1}^d {\rm ker} ( X^r - \l^r)^n \right).
\]
By Proposition \ref{prd1}, on a set of full measure each eigenvalue will have multiplicity one,
so we shall just talk about the spectrum.

Let $\tdn$ denote the commuting $d$-tuples of $n$-by-$n$ matrices that are upper triangular with respect to a fixed basis.
Any generic element $Q$ of $\tdn$ can be written uniquely as $A D A^{-1}$, where $D$ is a diagonal $d$-tuple and $A$ is an upper triangular $n$-by-$n$ matrix with $1$'s on the diagonal.
Let us write $\fno$ to denote the set of upper triangular $n$-by-$n$ matrices with $1$'s on the diagonal,
and ${\rm Ad}_A$ the linear operator on $\mn$ of conjugation by $A$; so ${\rm Ad}_A(X) = A X A^{-1}$.

\begin{lemma}
\label{lemd1}
Let $Q $ be in $\tdn \cap \comdng$.
Then the tangent space to $\vdn$ at $Q $ is 
\be
\label{eqd31}
\{ Z \in \mnd: [Q ^r, Z^s] = [Q ^s , Z^r], \quad \forall\ 1 \leq r \neq s \leq d \} .
\ee
The tangent space to $\tdn$ at $Q $ is $\{ Z \in \tdn : \eqref{eqd31}{\rm \ holds}\}$.
The real dimension of the first tangent space is $2n^2 + 2(d-1)n$, and of  the second is 
$n^2 + (2d-1)n$.
\end{lemma}
\bp
($\Rightarrow$) Suppose $\gamma(t)$ is a $C^1$ curve in $\vdn$ with $\gamma(0) = Q$
and $\gamma^\prime(0) = Z$.
Then $\gamma^r(t) = Q^r + t Z^r + o(t)$. As \[ 0 \= 
\gamma^r(t) \gamma^s(t) - \gamma^s(t) \gamma^r(t) \=
t\left(  [Q ^r, Z^s] - [Q ^s , Z^r]\right) + o(t) ,
\]
we conclude that $[Q ^r, Z^s] = [Q ^s , Z^r]$.

($\Leftarrow$) Suppose $Z$ is in \eqref{eqd31}.
Let $A$ be an upper triangular matrix that diagonalizes $Q$: $Q = A D A^{-1}$.
Let $Y = A^{-1}ZA$. Then
\be
\label{eqd51}
[ D^r, Y^s] \= [D^s, Y^r ] \quad \forall\ 1 \leq r \neq s \leq d .
\ee
Let $\l_j^r$ be the diagonal entries of $D^r$. Then \eqref{eqd51} says
\be
\label{eqd52}
Y^s_{ij} (\l^r_j - \l^r_i) \= Y^r_{ij} (\l^s_j - \l^s_i)\quad \forall\ 1 \leq r \neq s \leq d  .
\ee
Define a matrix $B$ by
\[
B_{ij} \= \begin{cases}
\frac{1}{\l^r_j - \l^r_i} Y^r_{ij}, \qquad &i \neq j
\\
0 \qquad& i = j \end{cases}.
\]
By \eqref{eqd52} the definition does not depend on $r$, and since $Q$ is generically diagonalizable
we never divide by zero. Let $E^r$ be the diagonal of $Y^r$.
Define a curve $\gamma(t)$ by
\[
\gamma^r(t) \= A e^{tB} (D^r + t E^r) e^{-tB} A^{-1} .
\]
Then $\gamma(t) \in \vdn$, and if $Z$ (and hence $Y$) is in $\tdn$, then so is $\gamma(t)$ for all $t$.
Moreover, $\gamma(0) = Q$ and 
\[
\gamma'(0) \= Ad_A ( [B,D] + E) \= Z.
\]
So the tangent space is as claimed.

To calculate the dimensions of the tangent spaces, note that $Ad_A$ is an isomorphism from
the tangent space at $D$ to the tangent space at $Q$, so it suffices to calculate the dimension of the tangent space at $D$. 
In $\vdn$, one can choose the diagonal entries freely, which gives $2dn$ dimensions, and 
the off diagonal entries for $r=1$, which by \eqref{eqd52} then determines all the other off-diagonal entries, for an additional $2(n^2-n)$ dimensions.
In the upper-triangular case, one can again choose the diagonal entries freely, and one has $n^2-n$ off diagonal dimensions.
\ep

Let us fix $ {Q_0} \in \tdn \cap \comdng$.
We shall use $T$ to mean tangent space.
Let $\vare$ be small and positive, and consider $\{ X : \| X - {Q_0} \|_F < \vare \}$.
For $\vare$ sufficiently small, the eigenvalues of $X$ will be so close to the eigenvalues of $Q_0$ that
 it makes sense to order the eigenvalues of $X$ so that the $j^{\rm th}$ eigenvalue is closest to the $j^{\rm th}$ eigenvalue of $Q_0$. Let $X = U Q U^*$ be an upper triangular factoring of $X$ 
 that preserves this order. The unitary 
 $U$ is unique up to multiplication by a diagonal unitary.
 
 As in Theorem \ref{thmc2} we can parametrize a neighborhood of $Q_0$ by
  \beq
 g: \tdn \times \exp (\azn) & \ \to \ & \vdn\\
 (Q, U) &\mapsto & U Q U^* \= {\rm Ad}_U  (Q) .
 \eeq
We want to decompose
\be
\label{eqd3}
\notag
dX \=
\kappa(Q) dQ d \nu
\ee
where $\nu$ is the probability
measure on the homogeneous space $\un/\T^n$ obtained by taking the quotient of Haar measure
on $\un$,
  and $dQ$ 
is Hausdorff measure of the appropriate dimension (which is $n^2 + (2d-1)n$) on $\tdn$.
For $d =1$ 
see the treatment by  Tao in \cite[Sec. 2.6.2]{ta12}.

Decompose $T_{Q_0} \vdn = T_{Q_0} \tdn \oplus (T_{Q_0} \tdn )^\perp.$
Let $P$ be orthogonal projection from $T_{Q_0} \vdn$ onto $T_{Q_0} \tdn$.
Let ${\rm ad}_{Q_0} (S)$ be the $d$-tuple $([Q_0^1, S], \dots, [Q_0^d,S])$.
Then we can write $dg$ in block form
\be
\label{eqd5}
dg_{(Q_0,1)}  \=
\bordermatrix{ & T_{Q_0} \tdn & \azn \cr
T_{Q_0} \tdn & {\rm id}  & - P\  {\rm ad}_{Q_0} \cr
(T_{Q_0} \tdn )^\perp & 0 & - P^\perp\ {\rm ad}_{Q_0}
 }
\ee
The Jacobian of $g$ is the square root of the absolute value of the determinant of the real linear map $(dg)^t dg$.
 From the structure of \eqref{eqd5}, this depends only on the $(2,2)$ entry.
Define a real linear map $\Lambda_Q : \azn \to \azn$ by
\be
\label{eqd6}
\notag
\la \Lambda_Q S_1, S_2 \ra \= \la P^\perp {\rm ad}_Q (S_1) , {\rm ad}_Q
 (S_2) \ra .
\ee
Then the Jacobian of $g$ at $(Q_0,1)$  is $\sqrt{|\det(\Lambda_{Q_0})|}$.
Summarizing we have proved the following.
\begin{proposition}
\label{prd17} Assume $\vdn$ is irreducible. With $X = U Q U^*$ as above, we have
\be
\label{eqd7}
\notag
w(X) dX \= \kappa(Q) w(Q) dQ d\nu ,
\ee
where $\kappa(Q) = \ccn \sqrt{|\det(\Lambda_Q)|}$.
\end{proposition}

\subsection{Calculating $\kappa$}
\label{subsecdc}

We can calculate $\kappa$ explicitly in the following special cases: $d=1$ (which was done 
by Ginibre \cite{gi65}); when $Q$ is diagonal; and when $n=2$.

Case: the $d$-tuple $Q$ equals $D$, a diagonal $d$-tuple. 
Let $\l_i$ denote the eigenvalues of $D$.
Now the tangent space $T_D \vdn$ decomposes into the direct sum of 
$T_D \tdn$ and the strictly lower
 triangular $d$-tuples $Z$ that satisfy \eqref{eqd31}.
 We have 
 \[
 {\rm ad}_D (E_{ji}) \= (\l_j^r - \l_i^r) E_{ji} .
 \]
 With $R_{ij}, S_{ij}$ as in the proof of Theorem \ref{thmc2}, we have that for $i < j$,
 \beq
 - P^\perp {\rm ad}_D : R_{ij} & \ \mapsto \ & - \frac{1}{\sqrt{2}} (\l_j^r - \l_i^r) E_{ji} \\
 S_{ij} & \ \mapsto \ &  \frac{i}{\sqrt{2}} (\l_j^r - \l_i^r) E_{ji} .
 \eeq
 Therefore the Jacobian of $g$ is 
 \be
 \label{eqd8}
 \notag
 \kappa(D) \= c_n \prod_{1 \leq i < j \leq n } |\l_i - \l_j|^2 ,
 \ee
 where $c_n = 2^{-n(n-1)/2} \ccn$.
 
 Case: $d=1$. When $d=1$ the conditions in \eqref{eqd31} are vacuous,
and $P^\perp$ is projection onto the strictly lower triangular matrices.
As shown in \cite[p. 190]{ta12}, if one orders the pairs $(j,i)$ with $ j > i$ in the order
$(n,1), (n,2), (n-1,1), (n,3), (n-1,2) ,(n-2,1), \dots$, then the $(2,2)$ entry of \eqref{eqd5} becomes
triangular, with diagonal entries ${Q}_{ii} - {Q}_{jj}$. Consequently 
\be
\label{eqd10}
\kappa(Q) \= c_n \prod_{1 \leq i < j \leq n } |\l_i - \l_j|^2 ,
\ee
where $\{ \l_i \}$ are the eigenvalues of $Q$.

Case: $n=2$. The space $T_{Q} \vdt \ominus T_{Q} \tdt$ is one dimensional.
If the eigenvalues of $Q$ are $\l_1$ and $\l_2$, then 
$T_{Q} \vdt \ominus T_{Q} \tdt$
is of complex dimension one, and is spanned by the $d$-tuple
\[
\omega \= 
\frac{1}{|\l_2 - \l_1|}
\left[
\begin{pmatrix}
0&0\\
(\l_2^r - \l_1^r) & 0 \\
\end{pmatrix}
\right] .
\]
Indeed,
\beq
[Q^s, \omega^r ] &\=&
\frac{1}{|\l_2 - \l_1|}
\begin{pmatrix}
0&0\\
(\l_2^r - \l_1^r) (\l_2^s - \l_1^s) & 0 \\
\end{pmatrix}\\
&=&
[Q^r, \omega^s],
\eeq
so by Lemma \ref{lemd1} $\omega$ is in $T_{Q} \vdt$. Moreover it is clearly orthogonal to
$T_{Q} \tdt$.

As 
\[
{\rm ad}_{Q_0} (E_{21}) \= \left[
\begin{pmatrix}
0&0\\
(\l_2^r - \l_1^r) & 0 \\
\end{pmatrix}
\right], 
\]
we get that
\beq
- P^\perp {\rm ad}_{Q_0}: R_{12} & \ \mapsto \ & \frac{1}{\sqrt{2}}| \l_2 - \l_1| \omega \\
S_{12} & \ \mapsto \ & \frac{-i}{\sqrt{2}} | \l_2 - \l_1 | \omega .
\eeq
Therefore \be
\label{eqd9}
\kappa(Q) \= c_2 |\l_2 - \l_1|^2 .
\ee

\subsection{From upper triangular to diagonal}

We are interested in determining the distribution $\rho_n(\l_1, \dots, \l_n)$ of the eigenvalues.
We shall mainly consider Gausssian weights on $\vdn$, that is weights
 of the form
\be
\label{eqd11}
w(X) \= C e^{- \gamma {\rm tr} \sum_{r=1}^d (X^r)^* X^r }\=  C e^{- \gamma \| X\|_F^2} .
\ee
The important property of Gaussian weights that we shall use, in addition to their unitary invariance, is that for upper triangular $d$-tuples $Q$  they factor as
\be
\label{eqd41}
w(Q) \= w_1 (\|D \|_F^2) w_2 ( \| Q - D \|_F^2) ,
\ee
where $D$ is the diagonal of $Q$.

When $d=1$, $\rho$ can be derived from $\kappa$ in a straightforward way.
Indeed, the off-diagonal terms can be chosen independently, and since the weight factors as in \eqref{eqd41}, 
 one can integrate out all the off-diagonal terms, and since $\kappa$ from \eqref{eqd10} depends only on the eigenvalues, we get that the density of the eigenvalues for $d=1$ is
\[
\rho_n(\l) \= 
C \prod_{1 \leq i < j \leq n } |\l_i - \l_j|^2 e^{ - \gamma \sum_{j=1}^n |\l_j|^2} .
\]
See \cite{gi65} or \cite[2.6.2]{ta12}.

For $d > 1$, we proceed as follows.
 Let $\diagdn$ denote the $d$-tuples of diagonal $n$-by-$n$ matrices.
 To make reading easier, we shall drop the subscript $0$ and fix an element $Q \in \tdn$.
 We have ${Q} = {\rm Ad}_A (D)$ for some $A \in \fno$ and  $D$ 
 the diagonal with the same eigenvalues as $Q$ (in the same order), which we 
 write $\l_1, \dots, \l_n$.
 We can parametrize a neighborhood of $Q$ in $\tdn$  by
  \beq
 G: \diagdn \times \fno  & \ \to \ & \tdn\\
 (\widetilde{D},\widetilde{A}) &\mapsto &  \widetilde{A} \widetilde{D} \widetilde{A}^{-1}  \=  {\rm Ad}_{\widetilde{A}} (\widetilde{D}) .
 \eeq

Let us calculate $dG$ at $(D, A)$. The tangent space of $\fno$ at $A$ is 
the space of strictly upper triangular matrices, which we shall denote by $\fnz$.
We shall work in complex coordinates, so the Jacobian will be $\det (dG)^* dG$
without a square root.
We have
\begin{eqnarray}
\nonumber 
dG_{(D,A)}: 
\cdn \times \fnz   & \to & T_{Q} \tdn \\
\nonumber
dG_{(D,A)} (E,B) &\= & \lim_{t \to 0} \frac{1}{t} [ (A + tB)(D + t E) (A + tB)^{-1}  - {Q} ] \\
&=&   {\rm Ad}_A( [A^{-1} B, D]) + {\rm Ad}_{A}(E) .
\label{eqd43}
\end{eqnarray}
As a basis for $\cdn$ we can choose $\{ \xi_p^r : 1 \leq p \leq n, 1 \leq r \leq d\}$,
where $\xi_p^r$ is the $d$-tuple that is $E_{pp}$ in the $r^{\rm th}$ slot and $0$ elsewhere.
As a basis for $\fnz$ we can choose $\{E_{ij} : 1 \leq i < j \leq n \} $.

We have $dG (\xi_p^r,0)$ is the $d$-tuple which in the $r^{\rm th}$ slot is the upper triangular $d$-tuple
that has $(k,l)$ entry $A_{kp} A^{-1}_{pl}$ in the $r^{\rm th}$ slot and $0$ elsewhere.
We have $dG(0, E_{ij})$ is the $d$-tuple whose $(k,l)$ entry in the  $r^{\rm th}$ slot is
\[
\delta_{ik} \l_j^r A^{-1}_{jl} - \sum_q A_{kq} \l_q^r A_{qi}^{-1} A^{-1}_{jl} .
\]
Define $\Gamma_Q$ to be the complex linear endomorphism of $\cdn \times \fnz$ given by
\be
\label{eqd12}
\Gamma_Q \= dG_{(D,A)}^* dG_{(D,A)} .
\ee
We have
\begin{eqnarray}
\nonumber
\lefteqn{\Gamma_Q ((\xi_{p_1}^{r_1}, E_{i_1 j_1}), (\xi_{p_2}^{r_2}, E_{i_2 j_2}))}\\
\label{eqd125}
\notag
&\=&
\sum_{k,l=1}^n \left[
\delta_{r_1, r_2} A_{k p_1} A^{-1}_{p_1 l} \overline{A_{k p_2}} \overline{A_{p_2 l}^{-1}}
\ + \ 
\sum_{q_2 =1}^n A_{k p_1} A^{-1}_{p_1 l} \overline{A_{k q_2}} \overline{\l^{r_1}_{q_2}}
\overline{A^{-1}_{q_2 i_2}} \overline{A^{-1}_{j_2 l}} \right.
\\
&+& \left. \sum_{q_1 =1}^n A_{kq_1} \l^{r_2}_{q_1} A^{-1}_{q_1 i_1} A^{-1}_{j_1 l} \overline{A_{kp_2}} \overline{A^{-1}_{p_2 l}} 
\ + \
\sum_{s=1}^d \sum_{q_1, q_2 =1}^n A_{kq_1} \l^{s}_{q_1} A^{-1}_{q_1 i_1} A^{-1}_{j_1 l}
\overline{A_{k q_2}} \overline{\l^{s}_{q_2}}
\overline{A^{-1}_{q_2 i_2}} \overline{A^{-1}_{j_2 l}}
\right] .
\nonumber
\end{eqnarray}
Then the Jacobian of $G$ is $| \det (\Gamma_Q)|$.
Thus we have proved the following (recall that the determinant of $\Lambda$ is calculated for a real linear map, and for $\Gamma$ it is complex linear, which explains why the first has a square root and the second does not).

\begin{theorem}
\label{thmd2} Assume $\vdn$ is irreducible. With $X = U A D A^{-1} U^*$ as above, we have
\be
\label{eqd13}
w(X) 
 dX \= \ccn  
 w_1 (\|D \|_F^2) w_2 ( \| Ad_A(D) - D \|_F^2) 
 \sqrt{|\det(\Lambda_{Ad_A(D)})|} \
 | \det ( \Gamma_{Ad_A(D)} ) | d\lambda dA d\nu .
 \ee
 The eigenvalue distribution of $X$ is given by
 \be
 \label{eqd135}
 \rho(\l) \= \ccn w_1 (\|D \|_F^2) \int_{\fno} w_2 ( \| Ad_A(D) - D \|_F^2) 
 \sqrt{|\det(\Lambda_{Ad_A(D)})|} \
 | \det ( \Gamma_{Ad_A(D)} ) |dA .
 \ee
\end{theorem}

\subsection{Density of eigenvalues}

We would like to find an explicit formula for the density of the eigenvalues, which would come from
 evaluating the right-hand side of \eqref{eqd135}. We shall do this first in the known case $d=1$, as an illustration. Then we shall do it in the case of $d$ arbitrary, but with $n=2$.

{\bf Case: $d=1$}. We have from \eqref{eqd43}
\beq
dG_{(D,A)} : \C^d \times \fnz & \ \to \ & T_Q \ton \\
(E,B) &\mapsto& Ad_A ( [A^{-1}B, D] + E) .
\eeq
Both the domain and codomain can be identified with the upper triangular matrices $\ton$.
As $A$ is upper triangular with $1$'s on the diagonal, multiplication by $A$ and $A^{-1}$ 
are unipotent linear operators on $\ton$, and hence have determinant $1$.
So the Jacobian of $G$ is the same as the Jacobian of the map
\[
(E,B) \ \mapsto \   [A^{-1}B, D] + E.
\]
Decomposing $\ton = \C^d \oplus \fnz$, this has the same determinant as
\[
B \ \mapsto  \ [A^{-1}B, D].
\]
As multiplication by $A^{-1}$ is unipotent on $\fnz$, this in turn has the same Jacobian as
\[
B \ \mapsto  \ [B, D], 
\]
which is 
\[
\Delta(\l) \ := \ \prod_{1 \leq i < j \leq n } |\l_i - \l_j|^2.
\]
So if we take \eqref{eqd13} and write it as
\[
w(X) dX \= w_1(\| D \|_F^2) w_2( \| A D A^{-1}  - D \|^2_F) \kappa(A D A^{-1}) 
\Delta(\l)
d\l d A d \nu,
\]
integrate with respect to $\nu$ and $A$ and use \eqref{eqd10}, we get
\be
\label{eqd435}
\rho(\l) \= c_n w_1(\| D \|^2_F) \Delta(\l) \int_{\fno} w_2( \| A D A^{-1}  - D \|^2_F) \Delta(\l) dA .
\ee
(We can pull $\kappa(A D A^{-1}) = c_n \Delta(\l)$ out of the integral since it does not depend on $A$.)
As 
\[
A D A^{-1}  - D \= [A,D] A^{-1} \= ad_A(D) A^{-1},
\]
 we can make the change of variables
 \beq
 B & \= & [A, D] A^{-1} \\
 dB &=& \Delta(\l) dA 
 \eeq
 to conclude that  the integral in 
 \eqref{eqd435} equals $\int_{\fnz} w_2(\| B\|^2_F) dB$, which is some constant.
 So we get Ginibre's density formula \cite{gi65} for some constant $C'$:
 \[
 \rho(\l) \= C' \Delta(\l) w_1(\| D\|_F^2)
 ,\]
 which if we assume $w$ is Gaussian as in  \eqref{eqd11} becomes
 \[
 \rho(\l) \= C' \Delta(\l)
  e^{-\gamma \sum_{j=1}^n |\l_j|^2} .
 \]
 
{\bf Case:} $n=2$, $d$ arbitrary. We shall restrict to the weight $w$ being as in \eqref{eqd11}. We can write
\[
A \= \begin{pmatrix}
1& \alpha \\
0& 1
\end{pmatrix} .
\]
From \eqref{eqd12}, the matrix for $\Gamma_Q$ has the form
\be
\label{eqd45}
\bordermatrix{~&\xi^r_1&\xi^r_2&E_{12}\cr
\xi^s_1& \delta_{rs} (1 + |\alpha|^2)  & -\delta_{rs} |\alpha|^2 & - \overline{\alpha}(\l_2^s - \l_1^s)
\cr
\xi^s_2 & -\delta_{rs} |\alpha|^2& \delta_{rs} (1 + |\alpha|^2) & \overline{\alpha}(\l_2^s - \l_1^s)
\cr
E_{12} &- \alpha (\bar \l_2^r - \bar \l_1^r) &\alpha (\bar \l_2^r - \bar \l_1^r) & |\l_2 - \l_1|^2
\cr}
\ee
To calculate the determinant of \eqref{eqd45} we shall use the Schur complement.
The first $2d$-by-$2d$ part of \eqref{eqd45} has determinant $(1 + 2|\alpha|^2)^d$.
The Schur complement with respect to the bottom right-hand entry is \[
|\l_2 - \l_1|^2 - \sum_{r=1}^d \la
\begin{pmatrix} 1+ |\alpha|^2 & - | \alpha|^2\\
-|\alpha|^2 & 1 + |\alpha|^2 \end{pmatrix}^{-1}
\begin{pmatrix}
 - \overline{\alpha}(\l_2^r - \l_1^r)\\
 \overline{\alpha}(\l_2^r - \l_1^r)
 \end{pmatrix} ,
 \begin{pmatrix}
 - \overline{\alpha}(\l_2^r - \l_1^r)\\
 \overline{\alpha}(\l_2^r - \l_1^r)
 \end{pmatrix}
\ra
\=
\frac{|\l_2 - \l_1|^2}{1 + 2 |\alpha|^2} .
\]
Multiplying these together, we get
\[
\det (\Gamma_Q) \= (1+ 2|\alpha|^2)^{d-1} |\l_2 - \l_1|^2 .
\]
Using this and \eqref{eqd9}, Theorem \ref{thmd2} gives
\[
w(X) dX \= C c_2\ e^{- \gamma(|\l_1|^2 + |\l_2|^2)} e^{- \gamma |\alpha|^2 |\l_2 - \l_1|^2}
|\l_2 - \l_1|^4 (1+ 2|\alpha|^2)^{d-1} d\lambda d \alpha d\nu .
\]
Integrate first with respect to $\nu$,then  use polar coordinates $\alpha = r e^{i\theta}$, and
 we get
\be
\label{eqd46}
\rho(\l) \= C' \ e^{- \gamma(|\l_1|^2 + |\l_2|^2)} |\l_2 - \l_1|^4 \int_0^\i e^{-\gamma |\l_2 - \l_1|^2 r^2}   (1 + 2r^2)^{d-1} r dr .
\ee
The following identity is elementary:
\be
\label{eqd47}
I_d \ := \ \int_0^\i r^{2d+1} e^{-kr^2} dr  \= \frac{d!}{2k^{d+1}} .
\ee
Using\eqref{eqd47} in \eqref{eqd46} with $k = \gamma |\l_2 - \l_1|^2$, we have proved the following result, where we absorb
some factors  depending on $d$  into $C'_d$.
\begin{theorem}
\label{thmd1}
Let $X$ be a random element of $\vdt$ chosen with distribution
\[
w(X) \= C e^{-\gamma \|X\|_F^2} .
\]
Then the eigenvalues $\l$  of $X$ in $(\C^d)^2$ have distribution
\be
\label{eqd48}
\rho^d(\l) \= 
C'_d \ e^{-\gamma ( |\l_1|^2 + |\l_2|^2)} |\l_2 - \l_1|^2 \sum_{j=0}^{d-1} \frac{1}{(d-1-j)!} \frac{2^j}{\gamma^{j+1}} \frac{1}{|\l_2 - \l_1|^{2j}} .
\ee
For $d=1,2,3$ we get respectively
\beq
\rho^1(\l) &\=& C_1^\prime e^{-\gamma ( |\l_1|^2 + |\l_2|^2)}\left[   \frac{1}{\gamma} |\l_2 - \l_1|^2\right]  \\
\rho^2(\l) &\=& C_2^\prime e^{-\gamma ( |\l_1|^2 + |\l_2|^2)} 
\left[ \frac{1}{\gamma} |\l_2 - \l_1|^2 + \frac{2}{\gamma^2}\right] \\
\rho^3(\l)  &\=& C_3^\prime e^{-\gamma ( |\l_1|^2 + |\l_2|^2)} 
\left[ \frac{1}{2\gamma} |\l_2 - \l_1|^2 + \frac{2}{\gamma^2}+ \frac{4}{\gamma^3} |\l_2 - \l_1|^{-2} \right] .
\eeq
\end{theorem}

\begin{remark}{\rm
The negative powers in \eqref{eqd48} are integrable, since $\l_j$ are in $\C^d$ and the smallest negative power is $4 - 2d$.
However, once $d \geq 3$, 
the negative powers show that, unlike in the Hermitian case, there is an {\em attraction} between the eigenvalues, not a repulsion.}
\end{remark}


%

\section{Notation}
\label{secnot}

For the convenience of the reader, we will list  our notation here.

\vs

$\acomdn \= \{ X \in \mnd: X^r X^s = - X^s X^r \ \forall\ r \neq s \}$

$\acomhdn \= \{ X \in \acomdn: (X^r)^* = X^r \ \forall \ r \}$

$
A_{n,\eta} \ =\ \{ y \in (\R^d)^n : \frac{1}{n^2} K_n(y) \leq E^Q + \eta \} 
$

 ${\rm Ad}_A(B) = A B A^{-1}$
 
  $\azn$  the skew-Hermitian matrices that have $0$'s on the diagonal
  
$\comdn \= \{ X \in \mnd: X^r X^s = X^s X^r \ \forall\ r,s \}$

$\comhdn \= \{ X \in \comdn: (X^r)^* = X^r \ \forall \ r \}$

$\comdng$ and $\comhdng$ denote the generically diagonalizable tuples (resp. the Hermitian ones)

${\rm diag}(A)$  the diagonal matrix whose entries are the  diagonal entries of $A$

$E_{ij}$ is the elementary matrix that has entry $1$ in the $(i,j)$ slot and $0$ elsewhere

$
\| X \|_F^2 \= \sum_{r=1}^d \sum_{i,j =1}^n |X_{ij}^r|^2 $ Frobenius norm

$\fnz$   upper triangular $n$-by-$n$ matrices with $0$'s on the diagonal

 $\fno$    upper triangular $n$-by-$n$ matrices with $1$'s on the diagonal

 $\mathcal G(d,n)$  the complex Grassmanian of all $d$-dimensional linear subspaces of $\C^n$

$GL(n,\C)$ is the invertible matrices

$
I^Q(\mu) \= \int \int \log \frac{1}{|x-y|} d\mu (x) d\mu(y) + \int Q(x) d \mu (x) 
$

$
k(s,t) \= \ \log \frac{1}{|s-t|} + \frac 12 ( Q(s) + Q(t) )
$

$
K_n(y) =  \sum_{1 \leq i < j \leq n} \log |y_i - y_j|^{-2} + (n-1) \sum_{i=1}^n Q(y_i) 
$

 $\mn$    $n$-by-$n$ matrices over $\C$
  
  $\mnd$ $d$-tuples in $\mn$ 
  
   $\niln \= \{ X \in \mn : ( X)^n = 0$
   
   $
\cp_n(y) \= D_n\,  e^{- n \sum_{j=1}^n Q(y_j)} \prod_{1 \leq i < j \leq n} |y_i - y_j|^2 
$
   
   $R_d$ are characteristic radii of a Gaussian distribution in $d$ variables (see Thm. \ref{emg})
   
 ${\mathfrak S}(n)$  the symmetric group on $n$ elements.

$\tdn$  the commuting $d$-tuples of $n$-by-$n$ matrices that are upper triangular with respect to a fixed basis
   
   ${\mathcal U}(n)$ is the $n$-by-$n$ unitary matrices
   
   $v(d)$ denotes the volume of the unit ball in $\R^d$


$\tw$ is the function of the eigenvalues $\l$ of $X \in \comhdn$ that satisfies $w(X) = \tw (\l)$

 $dX$ is  Haussdorf measure on ${\frak V}_n $

$\mu^Q$, 
the equilibrium measure for $Q$,  is the unique probability measure $\mu^Q$ that minimizes
$I^Q(\mu)$ (see Subsection \ref{subeq})

$\sigma^{d-1}_R$ is  normalized surface area measure
on the sphere $R S^{d-1}$

\bs
\bs
\bs
{\bf Acknowledgements.} This work was done while the author was visiting the Department of Mathematics at the University of California at Berkeley. The author would like to thank the Department for its hospitality.
The author would also like to thank the anonymous referee for many helpful suggestions that improved the readability of the paper.

\end{document}

%% file: def_latex5.tex
\def\mcc{M\raise.5ex\hbox{c}C}
\def\mccarthy{M\raise.5ex\hbox{c}Carthy}

\def\l{\lambda}

\def\vare{\varepsilon}

\let\i=\infty

\def\la{\langle}
\def\ra{\rangle}
\def\={\ = \ }

\def\C{\mathbb C}
\def\R{\mathbb R}
\def\T{\mathbb T}
\def\D{\mathbb D}

\def\be{\begin{equation}}
\def\ee{\end{equation}}
\def\beq{\begin{eqnarray*}}
\def\eeq{\end{eqnarray*}}
\def\vs{\vskip 5pt}
\def\bs{\vskip 12pt}
\def\bp{\begin{proof}}
\def\ep{\end{proof}}

\def\bl{\begin{lemma}}
\def\el{\end{lemma}}
\def\bt{\begin{theorem}}
\def\et{\end{theorem}}
\def\bprop{\begin{proposition}}
\def\eprop{\end{proposition}}
\def\bd{\begin{definition}}
\def\ed{\end{definition}}
\def\br{\begin{remark}}
\def\er{\end{remark}}
\def\bexer{\begin{exercise}}
\def\eexer{\end{exercise}}

\newtheorem{theorem}{Theorem}[section]
\newtheorem{proposition}[theorem]{Proposition}
\newtheorem{lemma}[theorem]{Lemma}
\newtheorem{corollary}[theorem]{Corollary}
\newtheorem{remark}[theorem]{Remark}

\newtheorem{definition}[theorem]{Definition}

\newtheorem{example}[theorem]{Example}
\newtheorem{question}[theorem]{Question}